\documentclass{amsart}

\usepackage{graphicx}
\usepackage{amsfonts}
\usepackage{amsmath}

\newtheorem{Thm}{Theorem}

\newtheorem{Lem}{Lemma}
\newtheorem{Prop}{Proposition}

\newtheorem{Def}{Definition}
\newtheorem{Rem}{Remark}

\newcommand{\z}{\mathbb{Z}}

\begin{document}
\title[Rank-two cluster algebras of affine type]{Combinatorial
interpretations for rank-two cluster algebras of affine type}
\author{Gregg Musiker}
\email{gmusiker@math.ucsd.edu}
\author{James Propp}
\email{James\_Propp@ignorethis.uml.edu} \thanks{{\rm 2000} {\it Mathematics
Subject Classification}. 05A99, 05C70}
\date{January 10, 2007}
\maketitle \tableofcontents

\begin{abstract}
Fomin and Zelevinsky \cite{ClustI} show that a certain
two-parameter family of rational recurrence relations, here called
the $(b,c)$ family, possesses the Laurentness property: for all
$b,c$, each term of the $(b,c)$ sequence can be expressed as a
Laurent polynomial in the two initial terms. In the case where the
positive integers $b,c$ satisfy $bc<4$, the recurrence is related to
the root systems of finite-dimensional rank $2$
Lie algebras; when $bc>4$, the recurrence is related to Kac-Moody
rank $2$
Lie algebras of general type \cite{Kac}.  Here we investigate the
borderline cases $bc=4$, corresponding to Kac-Moody Lie algebras of
affine type. In these cases, we show that the Laurent polynomials
arising from the recurence can be viewed as generating functions
that enumerate the perfect matchings of certain graphs. By providing
combinatorial interpretations of the individual coefficients of
these Laurent polynomials, we establish their positivity.
\end{abstract}

\section{Introduction}

In \cite{Laurent, ClustI}, Fomin and Zelevinsky prove that for all
positive integers $b$ and $c$, the sequence of rational functions
$x_n$ ($n \geq 0$) satisfying the ``$(b,c)$-recurrence''
$$x_n = \left\{ \begin{array}{ll}
(x_{n-1}^b+1)/x_{n-2} & \mbox{for $n$ odd} \\
(x_{n-1}^c+1)/x_{n-2} & \mbox{for $n$ even} \end{array} \right.$$ is
a sequence of Laurent polynomial in the variables $x_1$ and $x_2$;
that is, for all $n \geq 2$, $x_n$ can be written as a
sum of Laurent monomials of the form $a x_1^i x_2^j$, where the
coefficient $a$ is an integer and $i$ and $j$ are (not necessarily
positive) integers. In fact, Fomin and Zelevinsky conjecture that
the coefficients are always {\it positive\/} integers.




It is worth mentioning that variants of this recurrence typically
lead to rational functions that are not Laurent polynomials. For
instance, if one initializes with  $x_1,x_2$ and defines rational
functions
$$x_n = \left\{ \begin{array}{ll}
(x_{n-1}^b+1)/x_{n-2} & \mbox{for $n=3$ } \\
(x_{n-1}^c+1)/x_{n-2} & \mbox{for $n=4$ } \\
(x_{n-1}^d+1)/x_{n-2} & \mbox{for $n=5$ } \\
(x_{n-1}^e+1)/x_{n-2} & \mbox{for $n=6$ } \end{array} \right.$$ with
$b,c,d,e$ all integers larger than 1, then it appears that $x_5$ is
not a Laurent polynomial (in $x_1$ and $x_2$) unless $b=d$ and that
$x_6$ is not a Laurent polynomial unless $b=d$ and $c=e$. (This has
been checked by computer in the cases where $b,c,d,e$ are all
between 2 and 5.) 

One reason for studying $(b,c)$-recurrences is their relationship
with root systems associated to rank two
Kac-Moody Lie algebras.  Furthermore, algebras generated by a
sequence of elements satisfying a $(b,c)$-recurrence provide
examples of rank two cluster algebras, as defined in
\cite{ClustI,ClustII} by Fomin and Zelevinsky.
The property of being a sequence of Laurent polynomials,
\emph{Laurentness}, is in fact proven for all cluster algebras
\cite{ClustI} as well as a class of examples going beyond cluster
algebras \cite{Laurent}.  In this context, the positivity of the
coefficients is no mere curiosity, but is related to important
(albeit still conjectural) total-positivity properties of dual
canonical bases \cite{Canonical}.

The cases $bc<4$ correspond to finite-dimensional Lie algebras
(that is, semisimple Lie algebras), and these cases have been
treated in great detail by Fomin and Zelevinsky
\cite{ClustI,Sherman}. For example, the cases $(1,1)$, $(1,2)$, and
$(1,3)$ correspond respectively to the Lie algebras $A_2$, $B_2$,
and $G_2$. In these cases, the sequence of Laurent polynomials $x_n$
is periodic. More specifically, the sequence repeats with period 5
when $(b,c)=(1,1)$, with period 6 when $(b,c)=(1,2)$ or $(2,1)$, and
with period 8 when $(b,c)=(1,3)$ or $(3,1)$. For each of these
cases, one can check that each $x_n$ has positive integer
coefficients.

Very little is known about the cases $bc>4$, which should correspond
to Kac-Moody Lie algebras of general
type. It can be shown that for these cases, the sequence of Laurent
polynomials $x_n$ is non-periodic.

This article gives a combinatorial approach to the intermediate
cases $(2,2)$, $(1,4)$ and $(4,1)$, corresponding to Kac-Moody Lie
algebras of affine type; specifically algebras of types $A_1^{(1)}$
and $A_2^{(2)}$. Work of Sherman and Zelevinsky \cite{Sherman} has
also focused on the rank two affine case.  In fact, they are able to
prove positivity of the $(2,2)$-, $(1,4)$- and $(4,1)$-cases, as
well as a complete description of the positive cone. They prove both
cases simultaneously by utilizing a more general recurrence which
specializes to either case.  By using Newton polygons, root systems
and algebraic methods analogous to those used in the finite type
case \cite{ClustII}, they are able to construct the dual canonical
bases for these cluster algebras explicitly.

Our method is intended as a complement to the purely algebraic
method of Sherman and Zelevinsky \cite{Sherman}.
In each of the cases $(2,2)$, $(1,4)$ and $(4,1)$ we show that the
positivity conjecture of Fomin and Zelevinsky is true by providing
(and proving) a combinatorial interpretation of all the coefficients
of $x_n$. That is, we show that the coefficient of $x_1^i x_2^j$ in
$x_n$ is actually the cardinality of a certain set of combinatorial
objects, namely, the set of those perfect matchings of a particular
graph that contain a specified number of ``$x_1$-edges'' and a
specified number of ``$x_2$-edges''.  This combinatorial description
provides a different way of understanding the cluster variables, one
where the binomial exchange relations are visible geometrically.

The reader may already have guessed that the cases $(1,4)$ and
$(4,1)$ are closely related.  One way to think about this
relationship is to observe that the formulas
$$x_n = (x_{n-1}^b+1)/x_{n-2} \ \mbox{for $n$ odd}$$
and
$$x_n = (x_{n-1}^c+1)/x_{n-2} \ \mbox{for $n$ even}$$
can be re-written as
$$x_{n-2} = (x_{n-1}^b+1)/x_n \ \mbox{for $n$ odd}$$
and
$$x_{n-2} = (x_{n-1}^c+1)/x_n \ \mbox{for $n$ even};$$
these give us a canonical way of recursively defining rational
functions $x_n$ with $n<0$, and indeed, it is not hard to show
that
\begin{align} \label{recip} x_{-n}^{(b,c)}(x_1,x_2) =
x_{n+3}^{(c,b)}(x_2,x_1) \mathrm{~~~for~~~} n \in \z .\end{align} So
the $(4,1)$ sequence of Laurent polynomials can be obtained from the
$(1,4)$ sequence of Laurent polynomials by running the recurrence in
reverse and switching the roles of $x_1$ and $x_2$. Henceforth we
will not consider the $(4,1)$ recurrence; instead, we will study the
$(1,4)$ recurrence and examine $x_n$ for all integer values of $n$,
the negative together with the positive.

Our approach to the $(1,4)$ case will be the same as our approach to
the simpler $(2,2)$ case: in both cases, we will utilize perfect
matchings of graphs as studied in \cite[et
al.]{Aztec,Kuo,Reciprocity}.

\begin{Def} \rm For a graph $G = (V,E)$, which has an assignment of
weights $w(e)$ to its edges $e \in E$, a \emph{perfect matching} of
$G$ is a subset $S \subset E$ of the edges of $G$ such that each
vertex $v \in V$ belongs to exactly one edge in $S$.  We define the
\emph{weight} of a perfect matching $S$ to be the product of the
weights of its constituent edges, $$w(S) = \prod_{e\in S} w(e).$$
\end{Def}

With this definition in mind, the main result of this paper is the
construction of a family of graphs $\{G_n\}$ indexed by $n \in \z
\setminus \{1,2\}$ with weights on their edges such that the terms
of the $(2,2)$- (resp. $(1,4)$-) recurrence, $x_n$, satisfy $x_n =
p_n(x_1,x_2)/m_n(x_1,x_2)$; where $p_n(x_1,x_2)$ is the polynomial
$$\sum_{S \subset E \mathrm{~is~a~perfect~matching~of~}G_n} w(S),$$
and $m_n$ is the monomial $x_1^{c_1}x_2^{c_2}$ where $c_1$ and $c_2$
characterize the $2$-skeleton of $G_n$.  These constructions appear
as Theorem \ref{22thm} and Theorem \ref{14thm} in sections 2 and 3,
for the $(2,2)$- and $(1,4)$- cases, respectively.

%



Thus $p_n(x_1,x_2)$ may be considered a two-variable generating
function for the perfect matchings of $G_n$, and $x_n(x_1,x_2)$ may
be considered a generating function as well, with a slightly
different definition of the weight that includes ``global'' factors
(associated with the structure of $G_n$) as well as ``local''
factors (associated with the edges of a particular perfect
matching).  We note that the families $G_n$ with $n>0$ and $G_n$
with $n<0$ given in this paper are just one possible pair of
families of graphs with the property that the $x_n(x_1,x_2)$'s serve
as their generating functions.

The plan of this article is as follows. In Section 2, we treat the
case $(2,2)$; it is simpler than $(1,4)$, and makes a good warm-up.
In Section 3, we treat the case $(1,4)$ (which subsumes the case
$(4,1)$, since we allow $n$ to be negative). Section 4 gives
comments and open problems arising from this work.

\section{The $(2,2)$ case} \label{seq22}

Here we study the sequence of Laurent polynomials $x_1, x_2, x_3 =
(x_2^2+1)/x_1,$ etc.  If we let $x_1=x_2=1$, then the first few
terms of sequence $\{x_n(1,1)\}$ for $n\geq 3$ are
$2,5,13,34,89,\dots$.  It is not too hard to guess that this
sequence consists of every other Fibonacci number (and indeed this
fact follows readily from Lemma \ref{recur2} given on the next
page).

For all $n \geq 1$, let $H_n$ be the (edge-weighted) graph shown
below for the case $n=6$.

\bigskip

\begin{center}
\includegraphics[width = 3in , height = 0.9in]{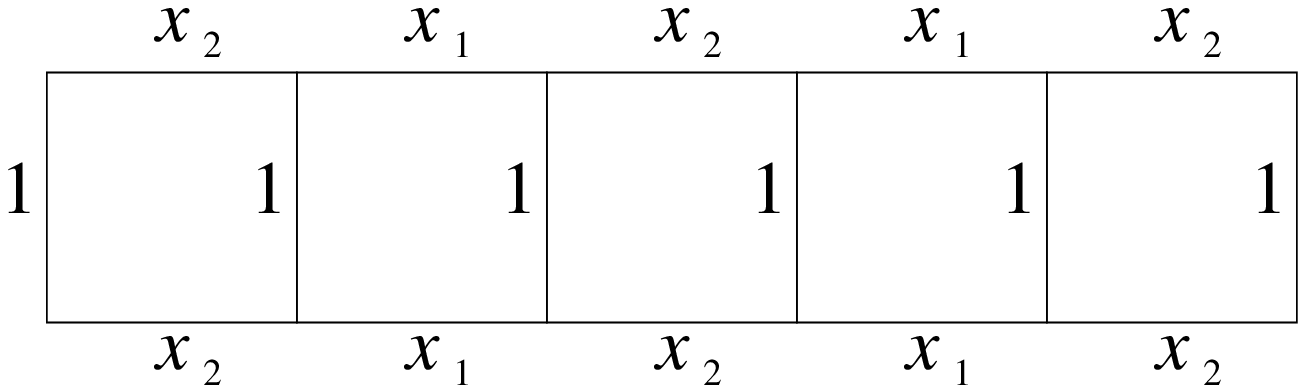} \\
The graph $G_5 = H_6$.
\end{center}

\bigskip

\noindent 
That is, $H_n$ is a 2-by-$n$ grid in which every vertical edge has
been assigned weight 1 and the horizontal edges alternate between
weight $x_2$ and weight $x_1$, with the two leftmost horizontal
edges having weight $x_2$, the two horizontal edges adjoining them
having weight $x_1$, the two horizontal edges adjoining {\it them\/}
having weight $x_2$, and so on (ending at the right with two edges
of weight $x_2$ when $n$ is even and with two edges of weight $x_1$
when $n$ is odd). Let $G_n = H_{2n-4}$ (so that for example the
above picture shows $G_5$), and let $p_n(x_1,x_2)$ be the sum of the
weights of all the perfect matchings of $G_n$. Also let
$m_n(x_1,x_2) = x_1^{n-2} x_2^{n-3}$ for $n \geq 3$.  We note the
following combinatorial interpretation of this monomial:
$m_n(x_1,x_2) = x_1^ix_2^j$ where $i$ is the number of square cells
of $G_n$ with horizontal edges having weight $x_2$ and $j$ is the
number of square cells with horizontal edges having weight $x_1$.
Using these definitions we obtain

\begin{Thm} \label{22thm}
For the case $(b,c)=(2,2)$, the Laurent polynomials $x_n$ satisfy
$$x_n(x_1,x_2)=p_n(x_1,x_2)/m_n(x_1,x_2) \mathrm{~for~} n \not = 1, 2$$
where $p_n$ and $m_n$ are given combinatorially as in the preceding
paragraph.
\end{Thm} E.g., for $n=3$, the graph $G_3 = H_2$ has two perfect
matchings with respective weights $x_2^2$ and $1$, so $x_3(x_1,x_2)
= (x_2^2+1)/x_1$.  For $n=4$, the graph $G_4=H_4$ has five perfect
matchings with respective weights $x_2^4$, $x_2^2$, $x_2^2$, $1$,
and $x_1^2$, so $p_4(x_1,x_2) = 1+2x_2^2+x_2^4+x_1^2$; since
$m_4(x_1,x_2) = x_1^2x_2$, we have $p_4(x_1,x_2)/m_4(x_1,x_2) =
(x_2^4+2x_2^2+1+x_1^2)/x_1^2x_2$, as required.

\begin{proof} We will have proved the claim if we can show that the Laurent
polynomials $p_n(x_1,x_2)/m_n(x_1,x_2)$ satisfy the same quadratic
recurrence as the Laurent polynomials $x_n(x_1,x_2)$; that is,
\begin{align} \label{ratexpr}
\frac{p_{n}(x_1,x_2)}{m_{n}(x_1,x_2)} \ {}
\frac{p_{n-2}(x_1,x_2)}{m_{n-2}(x_1,x_2)} =
\left(\frac{p_{n-1}(x_1,x_2)}{m_{n-1}(x_1,x_2)}\right)^2 +
1.\end{align}

\begin{Prop} \label{littleP}
The polynomials $p_n(x_1,x_2)$ satisfy the recurrence
\begin{align} \label{recur}
p_{n}(x_1,x_2) p_{n-2}(x_1,x_2) = (p_{n-1}(x_1,x_2))^2 + x_1^{2n-6}
x_2^{2n-8} \mathrm{~for~}n\geq 5.
\end{align}
\end{Prop}

\begin{proof}
To prove (\ref{recur}) we let $q_n(x_1,x_2)$ be the sum of the
weights of the perfect matchings of the graph $H_n$, so that
$p_n(x_1,x_2) = q_{2n-4}(x_1,x_2)$ for $n\geq 3$. Each perfect
matching of $H_n$ is either a perfect matching of $H_{n-1}$ with an
extra vertical edge at the right (of weight 1) or a perfect matching
of $H_{n-2}$ with two extra horizontal edges at the right (of weight
$x_1$ or weight $x_2$, according to whether $n$ is odd or even,
respectively). We thus have
\begin{eqnarray} q_{2n}   &=& q_{2n-1} + x_2^2 q_{2n-2} \\
q_{2n-1} &=& q_{2n-2} + x_1^2 q_{2n-3} \\ q_{2n-2} &=& q_{2n-3} +
x_2^2 q_{2n-4}.
\end{eqnarray} Solving the first and third equations for
$q_{2n-1}$ and $q_{2n-3}$, respectively, and substituting the
resulting expressions into the second equation, we get $(q_{2n} -
x_2^2 q_{2n-2}) = q_{2n-2} + x_1^2 (q_{2n-2} - x_2^2 q_{2n-4})$ or
$q_{2n} = (x_1^2+x_2^2+1) q_{2n-2} - x_1^2 x_2^2 q_{2n-4}$, so that
we obtain

\begin{Lem} \label{recur2}
\begin{align}
p_{n+1} = (x_1^2+x_2^2+1) p_n - x_1^2 x_2^2 p_{n-1}.
\end{align}
\end{Lem}

%

\noindent It is easy enough to verify that $$p_5p_3
=((x_2^2+1)^3+x_1^4+2x_1^2(x_2^2+1)\cdot (x_2^2+1)
 = \bigg((x_2^2+1)^2+x_1^2\bigg)^2 + x_1^4x_2^2 =  p_4^2 +
x_1^4x_2^2$$ so for induction we assume that
\begin{align} \label{induct}
p_{n-1}p_{n-3} = p_{n-2}^2 + x_1^{2n-8}x_2^{2n-10}
\mathrm{~for~}n\geq 5.
\end{align}
Using Lemma \ref{recur2} and (\ref{induct}) we are able to verify
that polynomials $p_n$ satisfy the quadratic recurrence relation
(\ref{recur}):
\begin{align*}
p_{n} p_{n-2} &= \\
(x_1^2+x_2^2+1)p_{n-1}p_{n-2} - x_1^2x_2^2p_{n-2}^2 &= \\
(x_1^2+x_2^2+1)p_{n-1}p_{n-2} - x_1^2x_2^2(p_{n-1}p_{n-3} - x_1^{2n-8}x_2^{2n-10}) &= \\
p_{n-1} ( (x_1^2 +x_2^2 +1)p_{n-2} -  x_1^2x_2^2p_{n-3}) +
x_1^{2n-6}x_2^{2n-8} &= p_{n-1}^2 + x_1^{2n-6}x_2^{2n-8}.
\end{align*}

\end{proof}



Since $m_n(x_1,x_2) = x_1^{n-2}x_2^{n-3}$ we have that recurrence
(\ref{ratexpr}) reduces to recurrence (\ref{recur}) of Proposition
\ref{littleP}.  Thus $p_n(x_1,x_2)/m_n(x_1,x_2)$ satisfy the same
initial conditions and recursion as the $x_n$'s, and we have proven
Theorem \ref{22thm} 
.
\end{proof}

An explicit formula has recently been found for the $x_n(x_1,x_2)$'s
by Caldero and Zelevinsky using the geometry of quiver
representations:

\begin{Thm} \cite[Theorem 4.1]{Caldero}, \cite[Theorem 2.2]{NewNote}
\label{explic}
\begin{align}
x_{-n} &=& \bigg(x_1^{2n+2} + \sum_{q+r \leq n} {n+1-r \choose
q}{n-q \choose r} x_1^{2q}x_2^{2r}\bigg)\bigg/x_1^nx_2^{n+1} \\
x_{n+3} &=& \bigg(x_2^{2n+2} + \sum_{q+r \leq n} {n-r \choose
q}{n+1-q \choose r} x_1^{2q}x_2^{2r}\bigg)\bigg/x_1^{n+1}x_2^{n}
\end{align} for all $n\geq 0$.
\end{Thm}
They also present expressions (Equations (5.16) of \cite{Caldero})
for the $x_n$'s in terms of Fibonacci polynomials, as defined in
\cite{Ysys}, which can easily seen to be equivalent to the
combinatorial interpretation of Theorem \ref{22thm}. Subsequently,
Zelevinsky has obtained a short elementary proof of these two
results \cite{NewNote}.

\subsection{Direct combinatorial proof of Theorem \ref{explic}}

Here we provide yet a third proof of Theorem \ref{explic}: instead
of using induction as in Zelevinsky's elementary proof, we use a
direct bijection.  This proof was found after Zelevinsky's result
came to our attention.  First we make precise the connection between
the combinatorial interpretation of \cite{Caldero} and our own.

\begin{Lem}
The number of ways to choose a perfect matching of $H_{m}$ with $2q$
horizontal edges labeled $x_2$ and $2r$ horizontal edges labeled
$x_1$ is the number of ways to choose a subset $S \subset
\{1,2,\dots, m-1\}$ such that $S$ contains $q$ odd elements, $r$
even elements, and no consecutive elements.
\end{Lem}

Notice that in the case $m=2n+2$, this number is the coefficient of
$x_1^{2r-n-1}x_2^{2q-n}$ in $x_{n+3}$ (for $n\geq 0$), and when
$m=2n+1$, this number is the coefficient of $x_1^{2r-n}x_2^{2q-n}$
in $s_n$, as defined in \cite{Caldero, Sherman, NewNote}.

\begin{proof}
There is a bijection between perfect matchings of $H_m$ and subsets
$S \subset \{1,2,\dots, m-1\}$, with no two elements consecutive.
We label the top row of edges of $H_m$ from $1$ to $m-1$ and map a
horizontal edge in the top row to the label of that edge.  Since
horizontal edges come in parallel pairs and span precisely two
vertices, we have an inverse map as well.
\end{proof}

With this formulation we now prove

\begin{Thm}
The number of ways to choose a subset $S \subset \{1,2,\dots,N\}$
such that $S$ contains $q$ odd elements, $r$ even elements,
and no consecutive elements is $${n+1-r \choose q}{n-q \choose r}$$
if $N=2n+1$ and $${n-r \choose q}{n-q \choose r}$$ if $N=2n$.

\end{Thm}

\begin{proof}

List the parts of $S$ in order of size and reduce the smallest by
$0$, the next smallest by $2$, the next smallest by $4$, and so on
(so that the largest number gets reduced by $2(q+r-1)$).

This will yield a multiset consisting of $q$ not necessarily
distinct odd numbers between $1$ and $2n+1-2(q+r-1) = 2(n-q-r+1)+1$
if $N=2n+1$, and between $1$ and $2(n-q-r+1)$ if $N=2n$, as well as
$r$ not necessarily distinct even numbers between $1$ and $2(n-q-r+1)$,
regardless of whether $N$ is $2n+1$ or $2n$.

Conversely, every such multiset, when you apply the bijection in
reverse, you get a set consisting of $q$ odd numbers and $r$ even
numbers in $\{1,2,\dots,2n\}$ (resp. $\{1,2,\dots, 2n+1\}$), no two
of which differ by less than $2$.

The number of such multisets is clearly $\left({n-q-r+1 \choose q}\right)
\times \left({n-q-r+1 \choose r}\right)$ in the first case
and $\left({n-q-r+2 \choose q} \right)\times \left({n-q-r+1
\choose r}\right)$ in the second case (since $2n+1$ is an
additional odd number), where $\left({n \choose k}\right) =$``$n$
multichoose $k$" $= {n+k-1 \choose k}$. Since $\left({n-q-r+1
\choose q} \right) \times \left({n-q-r+1 \choose r} \right) = {n-r
\choose q} {n-q \choose r}$ \bigg(resp. $\left({n-q-r+2 \choose q}
\right) \times \left({n-q-r+1 \choose r} \right) = {n+1-r \choose
q} {n-q \choose r}$\bigg), the claim follows.
\end{proof}

Once we know this interpretation for the coefficients of $x_{n+3}$
($s_n$), we obtain a proof of the formula for the entire sum, i.e.
Theorem \ref{explic} and Theorem 5.2 of \cite{Caldero} (Theorem 2.2 of \cite{NewNote}).

It is worth remarking that the extra terms $x_1^{2n+2}$ and $x_2^{2n+2}$
in Theorem \ref{explic} correspond to the extreme case in which one's
subset of $\{1,2,\dots,N\}$ consists of all the odd numbers in that range.


\subsection{Bijective proof of Lemma \ref{recur2}}
\label{bij}

The recurrence of Lemma \ref{recur2} can also be proven
bijectively by computing in two different ways the sum of the
weights of the perfect matchings of the graph $G_n \sqcup H_3$ (the
disjoint union of $G_n$ and $H_3$, which has all the vertices and
edges of graphs $G_n$ and $H_3$ and no identifications).  We provide
this proof since this method will be used later on in the $(1,4)$
case.


On the one hand, the sum of the weights of the perfect matchings of
$G_n$ is the polynomial $p_n$ and the sum of the weights of the
perfect matchings of $H_3$ is $x_1^2+x_2^2+1$, so the sum of the
weights of the perfect matchings of $G_n \sqcup H_3$ is
$(x_1^2+x_2^2+1)p_n$.

\bigskip

\begin{center}
\includegraphics[width = 4.5in , height = 0.9in]{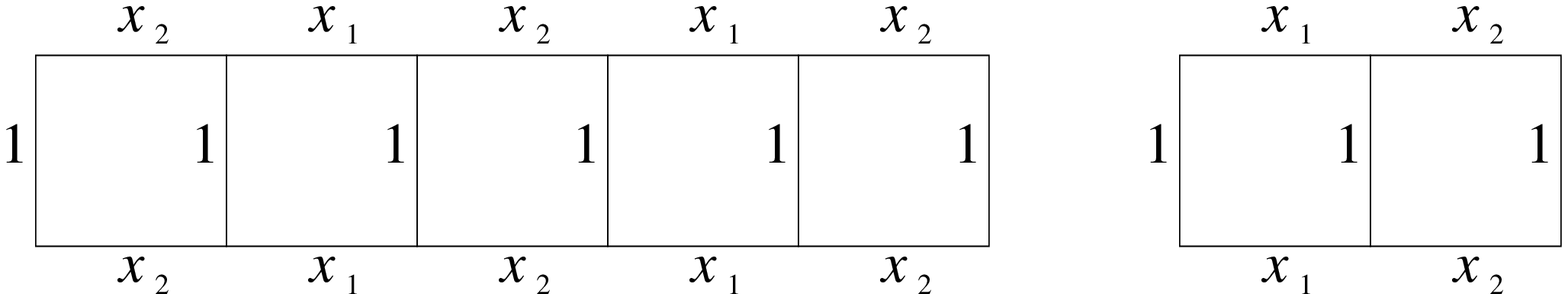} \\
The sum of the weights of all perfect matchings of $G_5 \sqcup H_3$
is $(x_1^2+x_2^2+1)p_5.$
\end{center}

\bigskip

On the other hand, observe that the graph $G_{n+1}$ can be obtained
from $G_n \sqcup H_3$ by identifying the rightmost vertical edge of
$G_n$ with the leftmost vertical edge of $H_3$. Furthermore, there
is a weight-preserving bijection $\phi$ between the set of perfect
matchings of the graph $G_{n+1}$ and the set of perfect matchings of
$G_n \sqcup H_3$ that do not simultaneously contain the two
rightmost horizontal edges of $G_n$ and the two leftmost horizontal
edges of $H_3$ (a set that
 can also be described as the set of
perfect matchings of $G_n \sqcup H_3$ that contain either the
rightmost vertical edge of $G_n$ or the leftmost vertical edge of
$H_3$ or both).  It is slightly easier to describe the inverse
bijection $\phi^{-1}$: given a perfect matching of $G_n \sqcup H_3$
that contains either the rightmost vertical edge of $G_n$ or the
leftmost vertical edge of $H_3$ or both, view the matching as a set
of edges and push it forward by the gluing map from $G_n \sqcup H_3$
to $G_{n+1}$.  We obtain a multiset of edges of $G_{n+1}$ that
contains either 1 or 2 copies of the third vertical edge from the
right, and then delete 1 copy of this edge, obtaining a \emph{set}
of edges that contains either 0 or 1 copies of that edge. It is not
hard to see that this set of edges is a perfect matching of
$G_{n+1}$, and that every perfect matching of $G_{n+1}$ arises from
this operation in a unique fashion. Furthermore, since the vertical
edge that we have deleted has weight 1, the operation is
weight-preserving.

The perfect matchings of $G_n \sqcup H_3$ that are not in the range
of the bijection $\phi$ are those that consist of a perfect matching
of $G_n$ that contains the two rightmost horizontal edges of $G_n$
and a perfect matching of $H_3$ that contains the two leftmost
horizontal edges of $H_3$. Removing these edges yields a perfect
matching of $G_{n-1}$ and a perfect matching of $H_1$.  Moreover,
every pair consisting of a perfect matching of $G_{n-1}$ and a
perfect matching of $H_1$ occurs in this fashion. Since the four
removed edges have weights that multiply to $x_1^2 x_2^2$, and $H_1$ has
just a single matching (of weight 1), we see that the perfect
matchings excluded from $\phi$ have total weight $x_1^2 x_2^2 p_{n-1}$  \cite[c.f.]{BenQuinn}.

\begin{Rem} \rm
We can also give a bijective proof of the quadratic recurrence
relation (\ref{recur}) by using a technqiue known as graphical
condensation which was developed by Eric Kuo \cite{Kuo}. He even
gives the unweighted version of this example in his write-up.
\end{Rem}

\begin{Rem} \rm As we showed via equation (\ref{recip}), there is a
reciprocity that allows us to relate the cluster algebras for the
$(b,c)$- and $(c,b)$-cases by running the recurrence backwards.
For the $(2,2)$-case, $b=c$ so we do not get anything new when we
run it backwards; we only switch the roles of $x_1$ and $x_2$.  This
reciprocity is a special case of the reciprocity that occurs not
just for $2$-by-$n$ grid graphs, but more generally in the problem
of enumerating (not necessarily perfect) matchings of $m$-by-$n$
grid graphs, as seen in \cite{Reciprocity2} and \cite{Reciprocity}.
For the $(1,4)$-case, we will also encounter a type of reciprocity.
\end{Rem}


\begin{Rem}  \rm
We have seen that the sequence of polynomials $q_n(x_1,x_2)$
satisfies the relation
$$q_{2n-4} q_{2n-8} = q_{2n-6}^2 + x_1^{2n-6} x_2^{2n-8}.$$
It is worth mentioning that the odd-indexed terms of the sequence
satisfy an analogous relation
$$q_{2n-3} q_{2n-7} = q_{2n-5}^2 - x_1^{2n-6} x_2^{2n-8}.$$
This relation can be proven via Theorem 2.3 of \cite{Kuo}.
%
%
%
%
%
%
%
%
%
In fact the sequence of Laurent polynomials $\{q_{2n+1}/x_1^nx_2^n:n
\geq 0\}$ are the collection of elements of the semicanoncial basis which
are not cluster monomials, i.e. not of the form $x_n^px_{n+1}^q$ for
$p,q \geq 0$. These are denoted as $s_n$ in \cite{Caldero} and
\cite{NewNote} and are defined as $\overline{S}_n(s_1)$ where
$\overline{S}_n(x)$ is the normalized Chebyshev polynomial of the
second kind, $S_n(x/2)$, and $s_1 = (x_1^2 + x_2^2 + 1)/ x_1x_2$. We
are thankful to Andrei Zelevinsky for alerting us to this fact.
We describe an analogous combinatorial interpretation for the $s_n$'s in the $(1,4)$-case in subsection \ref{14semi14}.
\end{Rem}

\newpage

\section{The $(1,4)$ case} \label{seq14}

In this case we let

\begin{align*} x_n &= { x_{n-1} + 1 \over x_{n-2}}
\mathrm{~~for~}n \mathrm{~odd}   \\
 &= {x_{n-1}^4 + 1 \over x_{n-2}} \mathrm{~~for~}n \mathrm{~even}
\end{align*}

\noindent for $n \geq 3$.  If we let $x_1 = x_2 = 1$, the first few
terms of $\{x_n(1,1)\}$
for $n \geq 3$ are: \\
\noindent $2, 17, 9, 386, 43, 8857, 206, 203321, 987, 4667522, 4729,
\dots$

\vspace{1em}

\noindent Splitting this sequence into two increasing subsequences,
we get for $n \geq 1$:

\begin{align} x_{2n+1} &= a_n = 2, 9, 43, 206, 987, 4729 \\
x_{2n+2} &= b_n = 17, 386, 8857, 203321, 4667522.
\end{align}

\vspace{1em}

\noindent Furthermore, we can run the recurrence backwards and
continue the sequence for negative values of $n$:

\noindent $\dots, 386, 9, 17, 2, 1, 1, 2, 3, 41, 14, 937, 67, 21506,
321, 493697, 1538, 11333521, 7369 \dots$ \\whose negative terms
split into two increasing subsequences (for $n \geq 1$)

\begin{align} x_{-2n+2} &= c_n =  2, 41, 937, 21506, 493697, 11333521 \\
x_{-2n+1} &= d_n = 3, 14, 67, 321, 1538, 7369.
\end{align}

\vspace{1em} \noindent As in the $(2,2)$-case, it turns out that
this sequence $\{x_n(1,1)\}$ (respectively $\{x_n\}$) has a
combinatorial interpretation as the number (sum of the weights) of
perfect matchings in a sequence of graphs.  We prove that these
graphs, which we again denote as $G_n$, have the $x_n$'s as their
generating functions in the later subsections.  We first give the
unweighted version of these graphs where graph $G_n$ contains
$x_n(1,1)$ perfect matchings.  We describe how to assign weights to
yield the appropriate Laurent polynomials $x_n$ in the next
subsection, deferring proof of correctness until the ensuing two
subsections. The proof of two recurrences, in sections
\ref{sectfirst} and \ref{sectsecond}, will conclude the proof of
Theorem \ref{14thm}. The final subsection provides a combinatorial
interpretation for elements of
the semicanonical basis 
that
are distinct from cluster monomials.

\begin{Def} \label{14graphs} \rm We will have four types of graphs $G_n$, one for each of the
above four sequences (i.e.\ for $a_n$, $b_n$, $c_n$, and $d_n$).
Graphs in all four families are built up from squares (consisting of
two horizontal and two vertical edges) and octagons (consisting of
two horizontal, two vertical, and four diagonal edges), along with
some extra arcs.  We describe each family of graphs by type.

Firstly, $G_3$ ($a_1$) is a single square, and $G_5$ ($a_2$) is an
octagon surrounded by three squares.  While the orientation of this
graph will not affect the number of perfect matchings, for
convenience of describing the rest of the sequence $G_{2n+3}$, we
assume the three squares of $G_5$ are attached along the eastern,
southern, and western edges of the octagon and identify $G_3$ with
the eastern square.  For $n \geq 3$ the graph associated to
$a_{n+1}$, $G_{2n+3}$, can be inductively built from the graph for
$a_n$, $G_{2n+1}$ by attaching a complex consisting of one octagon
with two squares attached at its western edge and northern/southern
edge (depending on parity).  We attach this complex to the western
edge of $G_{2n+1}$, and additionally adjoin one arc between the
northeast (resp. southeast) corner of the southern (resp. northern)
square of $G_{2n+3} \setminus G_{2n+1}$ and the southeast (resp.
northeast) corner of the northern (resp. southern) square of
$G_{2n+1} \setminus G_{2n-1}$.

We can inductively build up the sequence of graphs corresponding to
the $d_n$'s, $G_{-2n-1}$, analogously.  Here $G_{-1}$, consists of a
single octagon with a single square attached along its northern
edge.  We attach the same complex (one octagon and two squares)
except this time we orient it so that the squares are along the
northern/southern and \emph{eastern} edges of the octagon.  We then
attach the complex so that the eastern square attaches to
$G_{-2n+1}$.  Lastly we adjoin one arc between the north\emph{west}
(resp. south\emph{west}) corner of the southern (resp. northern)
square of $G_{-2n-1} \setminus G_{-2n+1}$ and the south\emph{west}
(resp. north\emph{west}) corner of the northern (resp. southern)
square of $G_{-2n+1} \setminus G_{-2n+3}$.

The graph corresponding to $b_1$, $G_{4}$, is one octagon surrounded
by four squares while the graph corresponding to $c_1$, $G_0$,
consists of a single octagon. As in the case of $a_n$ or $d_n$, for
$n \geq 1$, the graphs for $b_{n+1}$ ($G_{2n+4}$) and $c_{n+1}$
($G_{-2n}$) are constructed from $b_n$ and $c_n$ (resp.), but this
time we add a complex of an octagon,
two squares and an arc on \emph{both sides}.  Note that this gives
these graphs symmetry with respect to rotation by $180^\circ$.

The graphs $G_{2n+2} ~(b_n)$ consist of a structure of octagons and
squares such that there are squares on the two ends, four squares
around the center octagon, and additional arcs shifted towards the
center (the vertices joined by an arc lie on the vertical edges
closest to the central octagon). On the other hand, the graphs
$G_{-2n+2} ~(c_n)$ have a structure of octagons and squares such
that the central and two end octagons have only two squares
surrounding them, and the additional arcs are shifted towards the
outside.
\end{Def}

For the reader's convenience we illustrate these graphs for small
$n$.  In our pictures, the extra arcs are curved, and the other edges 
are line segments.  It should be noted that all these graphs are
planar, even though it is more convenient to draw them in such a way
that the (curved) extra arcs cross the other (straight) edges.

\newpage \setlength{\unitlength}{1in}

\begin{tabular}{rrrr}
$n$ & $x_n(1,1)$ & Type &  \hspace{1.3in} Graph $~~~~~~~G_n$ \\
\hline \vspace{0.1em} \\

$-8$ & $493697 $ & $c_5$ & \put(-1.8,-.25){\includegraphics[width =
3in ,
height = .5in]{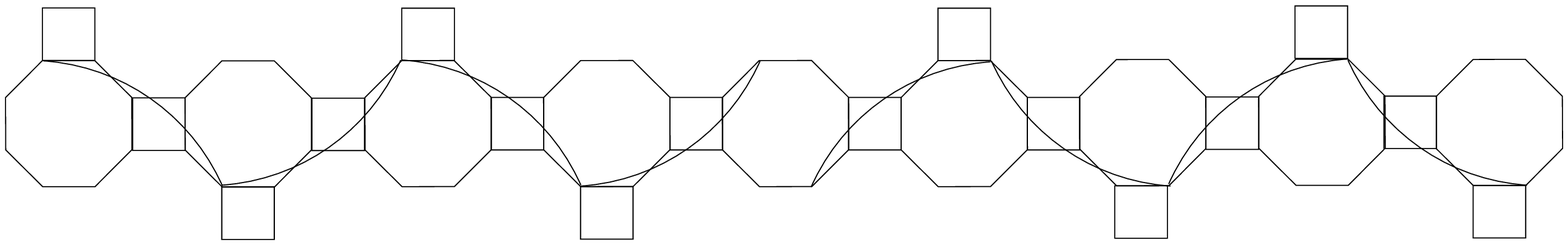}}  \\

\vspace{0.2em} \\

$-7$ & $321 $ & $d_4$ & \put(-1.05,-.25){\includegraphics[width =
1.5in ,
height = .5in]{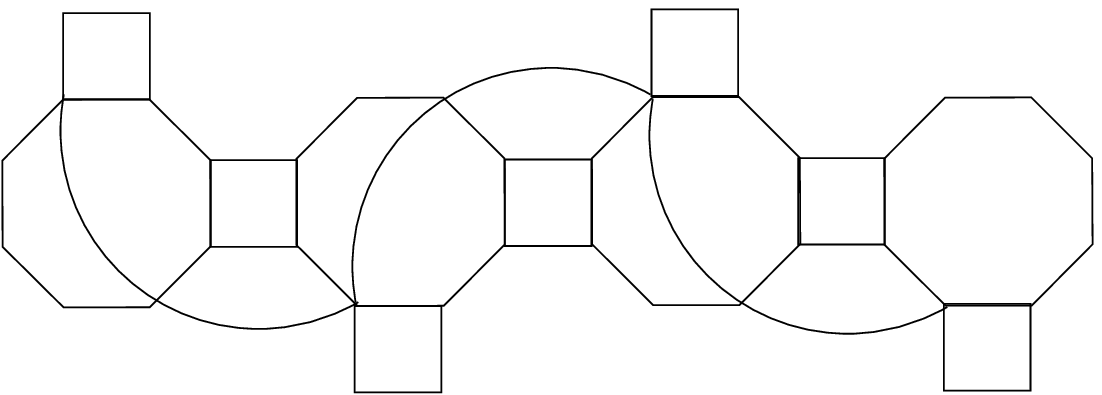}} \\

\vspace{0.2em} \\

$-6$ & $21506 $ & $c_4$ & \put(-1.55,-.25){\includegraphics[width =
2.5in ,
height = .5in]{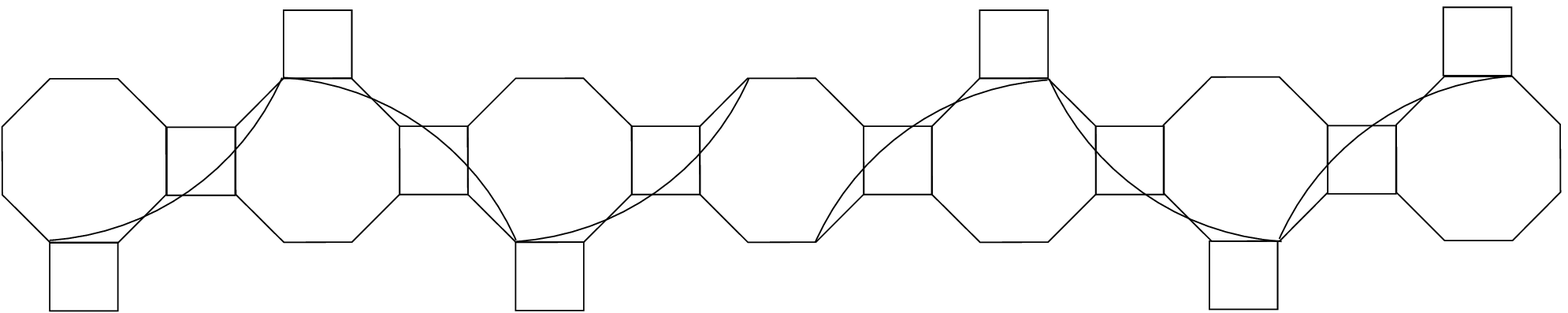}} \\

\vspace{0.2em} \\

$-5$ & $67 $ & $d_3$ &
 \put(-.8,-.25){\includegraphics[width = 1in , height = .5in]{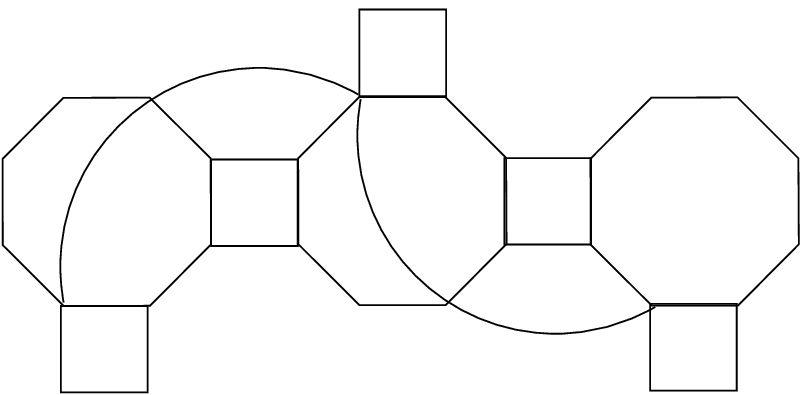}} \\

\vspace{0.2em} \\

$-4$ & $937 $ & $c_3$ &
\put(-1.3,-.25){\includegraphics[width = 2in , height = .5in]{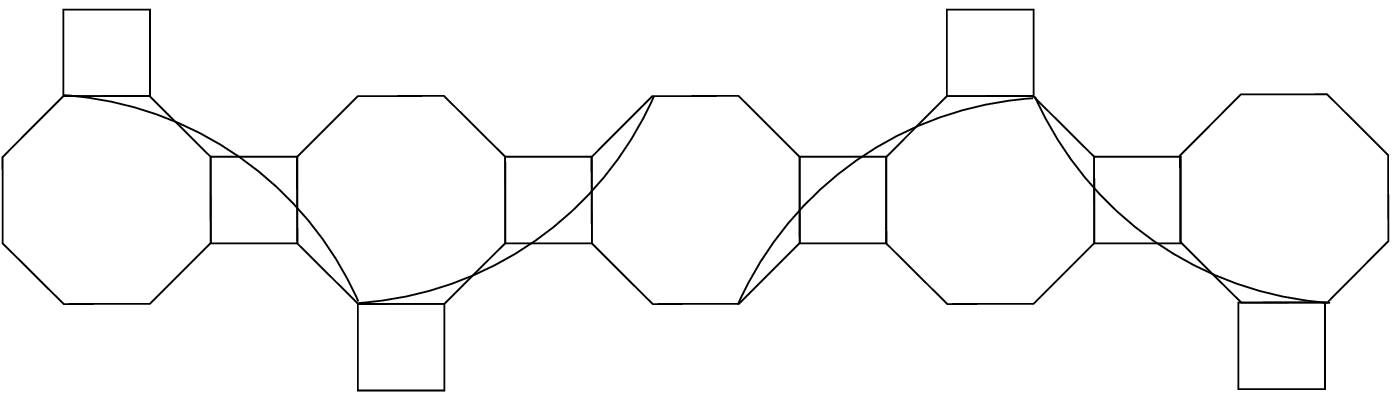}} \\

\vspace{0.2em} \\

$-3$ & $14 $ & $d_2$ &
\put(-.675,-.25){\includegraphics[width = 0.75in , height = .5in]{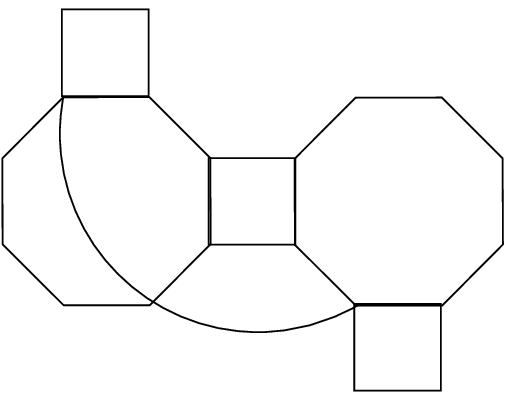}} \\

\vspace{0.2em} \\

$-2$ & $41 $ & $c_2$ &
\put(-.8,-.25){\includegraphics[width = 1in , height = .5in]{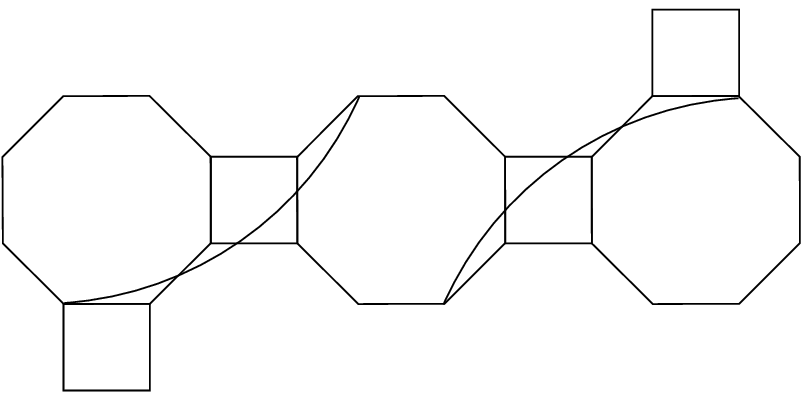}} \\

\vspace{0.2em} \\

$-1$ & $3 $ & $d_1$ &
\put(-.5,-.25){\includegraphics[width = .4in , height = .5in]{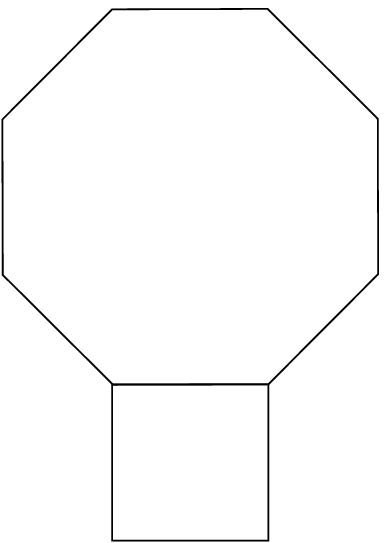}} \\

\vspace{0.2em} \\

$0$ & $2$ & $c_1$ &
\put(-.55,-.25){\includegraphics[width = .5in , height = .5in]{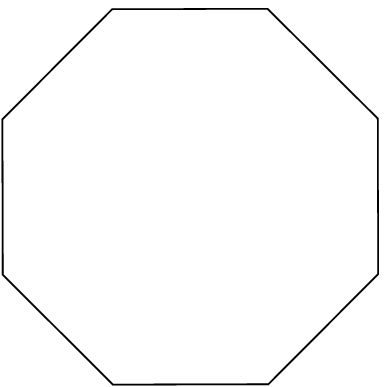}} \\

\vspace{0.2em} \\

$1$ & $1 $ & &
\put(-.4,0){$x_1$} \\

\vspace{0.2em} \\

$2$ & $1 $ & &
\put(-.4,0){$x_2$} \\
\end{tabular}

\newpage
\begin{tabular}{rrrr}
\centering $3$ & $2 $ & $a_1$ &
\put(1.875,-.075){\includegraphics[width = .25in , height = .25in]{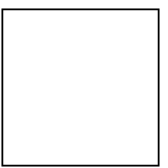}} \\

\vspace{0.2em} \\

$4$ & $17 $ & $b_1$ &
\put(1.75,-.2){\includegraphics[width = .5in , height = .5in]{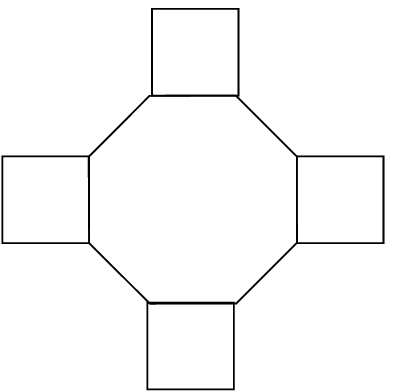}} \\

\vspace{0.2em} \\

$5$ & $9 $ & $a_2$ &
\put(1.75,-.2){\includegraphics[width = .5in , height = .5in]{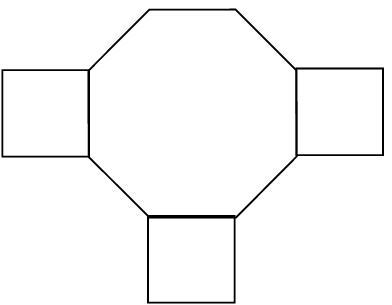}} \\

\vspace{0.2em} \\

$6$ & $386 $ & $b_2$ &
\put(1.4,-.2){\includegraphics[width = 1.2in , height = .5in]{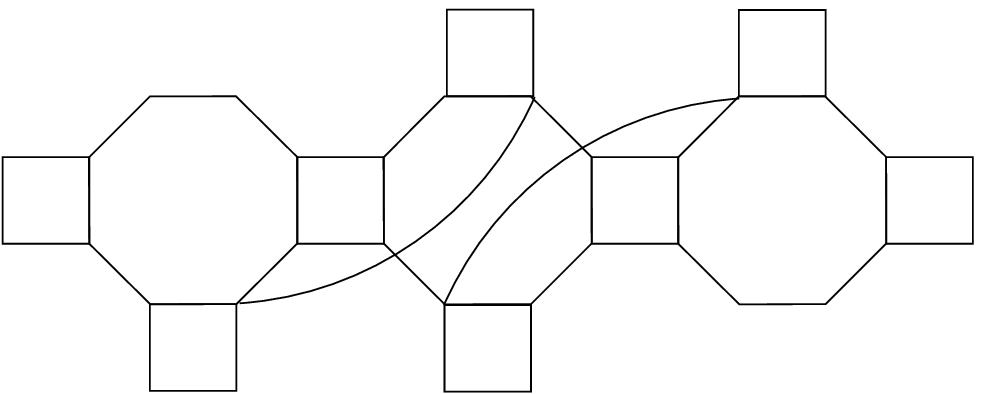}} \\

\vspace{0.2em} \\

$7$ & $43 $ & $a_3$ &
\put(1.6,-.2){\includegraphics[width = 0.8in , height = .5in]{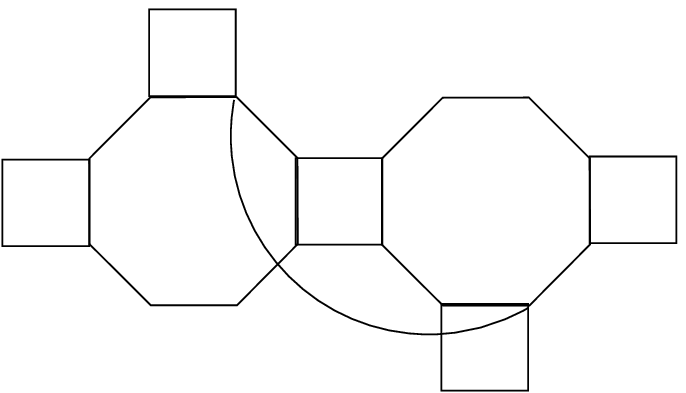}} \\

\vspace{0.2em} \\

$8$ & $8857 $ & $b_3$ &
\put(1.1,-.2){\includegraphics[width = 1.8in , height = .5in]{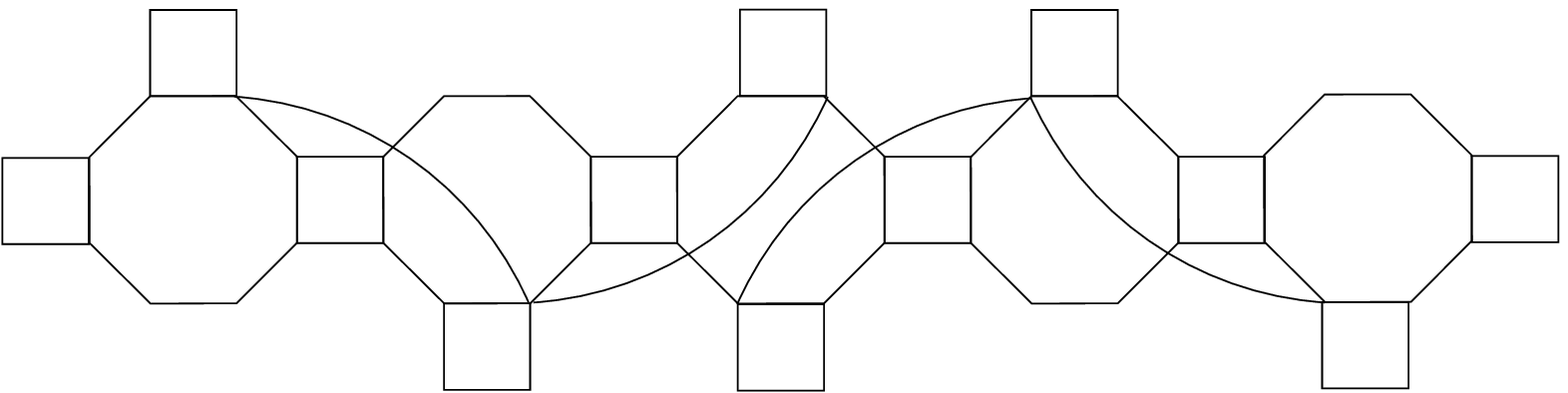}} \\

\vspace{0.2em} \\

$9$ & $206 $ & $a_4$ &
\put(1.4,-.2){\includegraphics[width = 1.2in , height = .5in]{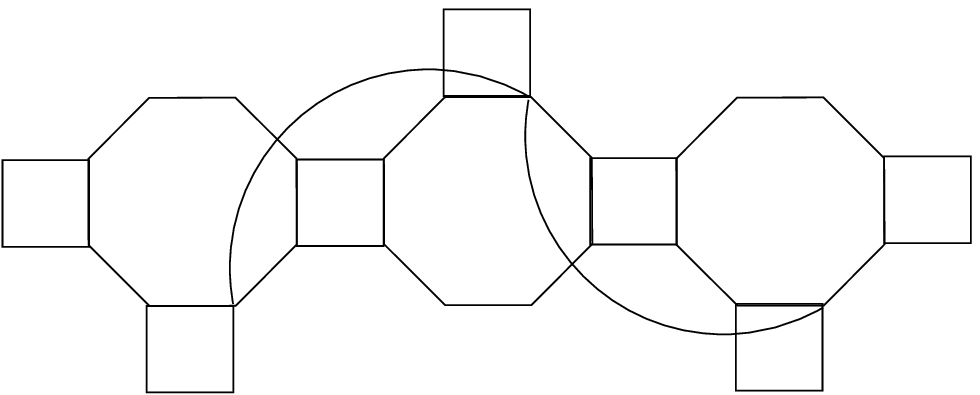}} \\

\vspace{0.2em} \\

$10$ & $203321 $ & $b_4$ &
\put(.6,-.2){\includegraphics[width = 2.8in , height = .5in]{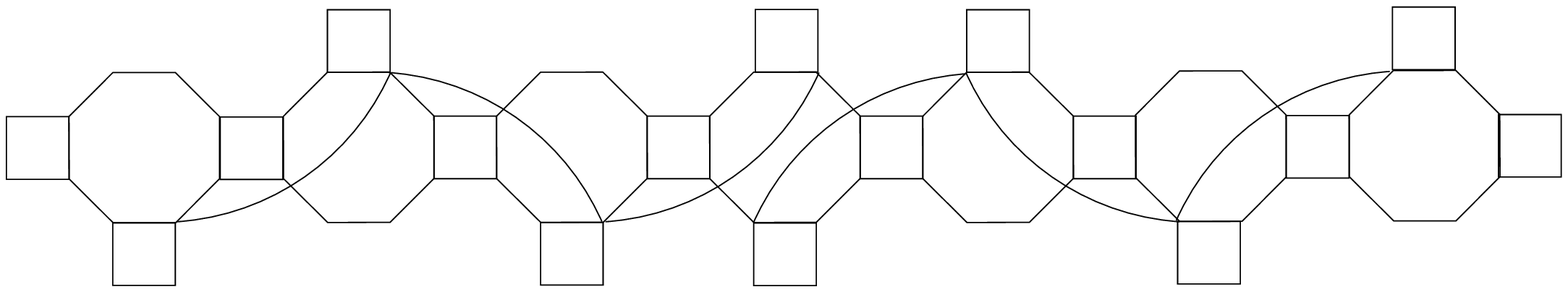}} \\

\vspace{0.2em} \\

$11$ & $987 $ & $a_5$ &
\put(1.15,-.2){\includegraphics[width = 1.7in , height = .5in]{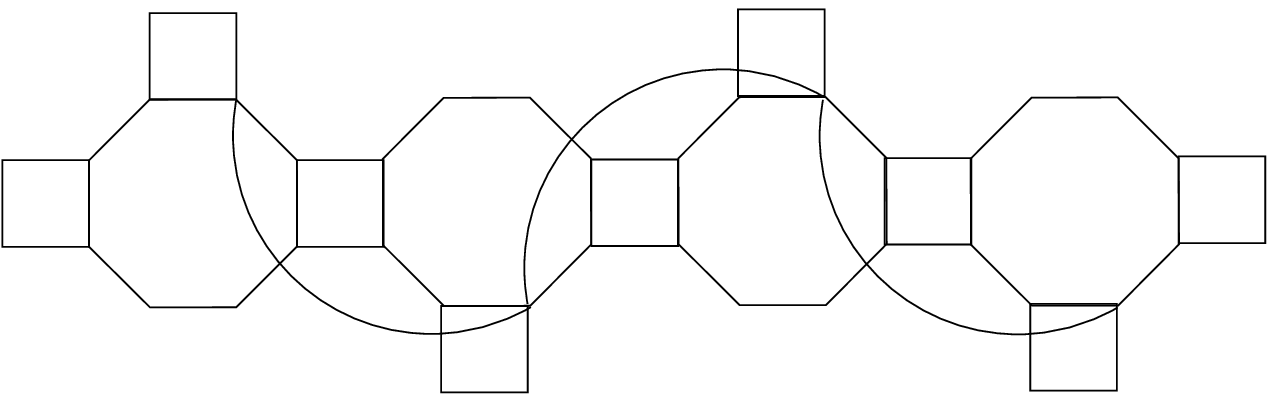}} \\

\vspace{0.2em} \\

$12$ & $4667522$ & $b_5$ &
\put(.4,-.2){\includegraphics[width = 3.2in, height = .5in]{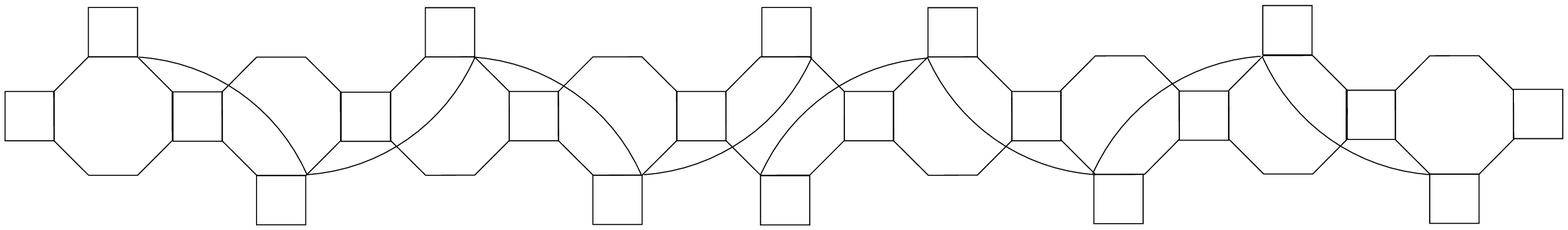}} \\

\vspace{0.2em} \\

$13$ & $4729 $ & $a_6$ &
\put(.9,-.2){\includegraphics[width = 2.2in , height = .5in]{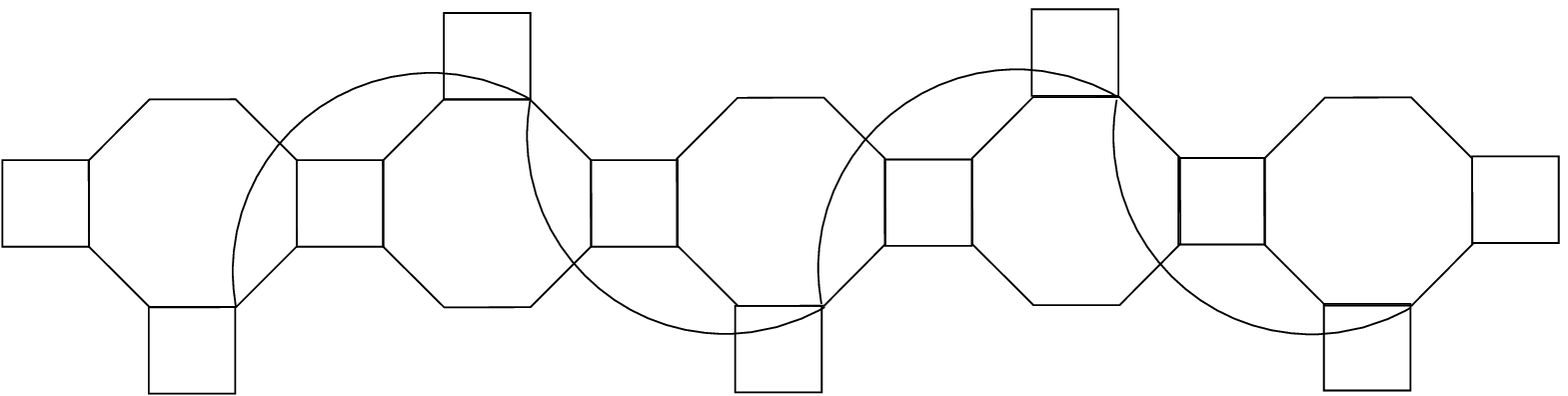}} \\

\end{tabular}

\newpage

\begin{Rem} \rm
As described above, the sequence of graphs corresponding to the
$a_n$'s and the $d_n$'s can both be built up inductively.  In fact,
if one assumes that the graphs associated with $d_n$ are ``negative"
then one can even construct $a_1$ from $d_1$ by ``adding" two
squares and an octagon.  The negative square and octagon cancels
with the positive square and octagon, leaving only a square for the
graph of $a_1$.  (When we construct the graph associated to $a_2$
from the graph for $a_1$, we do not add an arc as in the $n \geq 2$
case. Similarly, we omit an extra arc when we construct the graph
for $a_1$ from the graph for $d_1$. We do not have a ``principled''
explanation for these exceptions, but we do note that in these two
cases, the graph is sufficiently small that there are no candidate
vertices to connect by such an arc.)

Comparing graphs with equal numbers of octagons, we find a nice
reciprocity between the graph for $x_{2n+3}$ and the graph for
$x_{-2n+1}$ for $n \geq 1$. Namely, the two graphs are isomorphic
(up to horizontal and vertical reflection) except for the fact that
the graph $G_{2n+3}$ contains two squares on the left and right ends
while graph $G_{-2n+1}$ lacks these squares. Notice that graph $G_3$
lies outside this pattern since it contains no octagons and instead
contains a single square.

Studying the other types of graphs, i.e. those corresponding to
$b_n$ and $c_n$, we see that in each pair of graphs containing the
same number of octagons, there is a nice reciprocity between the
two. (Compare the definitions of the graphs associated to the
$b_n$'s with those associated to the $c_n$'s.)  These two reciprocal
relationships reduce into one, i.e.\ $G_{-n}$ and $G_{n+4}$ are
reciprocal graphs for all integers $n \geq 0$.  Recall that a
reciprocity also exists (between graphs $G_{-n}$ and $G_{n+3}$) for
the (2,2) case, so perhaps these reciprocities are signs of a more
general phenomenon. \end{Rem}

\subsection{Weighted versions of the graphs}

\label{weighted}

We now turn to the analysis of the sequence of Laurent polynomials
$x_n(x_1,x_2)$, and give the graphs $G_n$ weights on the various
edges. 
In this case, the denominator depends on the number of faces in the
graph, ignoring extra arcs. The exponent of $x_1$ in the denominator
will equal the number of squares while the exponent of $x_2$ will
equal the
number of octagons. Because of this interpretation, for $n \not = 1$
or $2$ we will rewrite $x_n$ as ${p_n(x_1,x_2) \over
x_1^{sq(n)}x_2^{oct(n)}}$ where $sq(n)$ and $oct(n)$ are both
nonnegative integers.  By the description of graphs $G_n$, we find

\begin{eqnarray} \label{sqn} sq(n) &=& \left\{ \begin{array}{ll}
|n-1|-1 & \mbox{for $n$ odd} \\
|2n-2|-2 & \mbox{for $n$ even} \end{array} \right. \\
\label{octn} oct(n) &=& \left\{ \begin{array}{ll}
|{n \over 2}-1|-{1 \over 2} & \mbox{for $n$ odd} \\
|n-2|-1 & \mbox{for $n$ even}. \end{array} \right.
\end{eqnarray}

To construct these weighted graphs, we take the graphs $G_n$ and
assign weights such that each of the squares has one edge of
weight $x_2$ and three edges of weight $1$ while the octagons have
weights alternating between $x_1$ and $1$.  As an example, consider
the following close-up of the graph associated to $x_{10}$.
\newpage

\begin{center} \includegraphics[width = 5in , height = 1in]{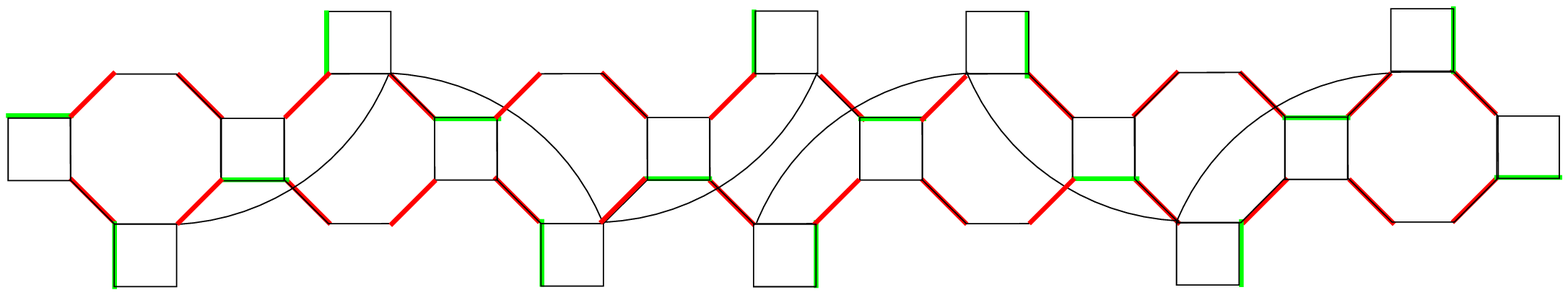} \\
The graph associated to $x_{10}$. \end{center}

\noindent The vertical and horizontal edges colored in green are
given weight $x_2$, and the diagonal edges marked in red are given
weight $x_1$.  All other edges are given weight $1$. Notice that the
vertical edges weighted $x_2$ lie furthest away from the arcs and the
horizontal edges weighted $x_2$ alternate between top and bottom,
starting with bottom on the righthand side.

\vspace{2em} Table of $x_n$ for small $n$:

\vspace{2em}
\begin{tabular}{cc}
$n$ & $x_n$
\\ \hline
\vspace{0.1em} \\
$-3$ & ${(x_2+1)^3 + 2x_1^4+ 3x_1^4x_2+x_1^8  \over x_1^3x_2^2}$ \\
\vspace{0.1em} \\
$-2$ & ${(x_2+1)^4 + 3x_1^4 +8x_1^4x_2+6x_1^4x_2^2 +  3x_1^8+4x_1^8x_2 +  x_1^{12} \over x_1^4x_2^3}$ \\
\vspace{0.1em} \\
$-1$ & ${(x_2+1)+x_1^4 \over x_1x_2}$ \\
\vspace{0.1em} \\
$0$ & ${x_1^4+1 \over x_2}$ \\
\vspace{0.1em} \\
$1$ & $x_1$ \\
\vspace{0.1em} \\
$2$ & $x_2$ \\
\vspace{0.1em} \\
$3$ & ${x_2+1 \over x_1}$ \\
\vspace{0.1em} \\
$4$ & ${(x_2+1)^4+x_1^4 \over x_1^4x_2}$ \\
\vspace{0.1em} \\
$5$ & ${(x_2+1)^3+x_1^4 \over x_1^3x_2}$ \\
\vspace{0.1em} \\
$6$ & ${(x_2+1)^8 + 3x_1^4 + 16x_1^4x_2 + 34x_1^4x_2^2 +
36x_1^4x_2^3 + 19x_1^4x_2^4
+ 4x_1^4x_2^5 + 3x_1^8 + 8x_1^8x_2 + 6x_1^8x_2^2 + x_1^{12} \over x_1^8x_2^3}$ \\
\vspace{0.1em} \\
$7$ & ${(x_2+1)^5 + 2x_1^4 + 5x_1^4x_2 + 3x_1^4x_2^2 + x_1^8 \over
x_1^5x_2^2}$
\end{tabular}\vspace{2em}





\noindent We now wish to prove the following.

\begin{Thm} \label{14thm}
For the case $(b,c)=(1,4)$, the Laurent polynomials $x_n$ satisfy
$$x_n(x_1,x_2)=p_n(x_1,x_2)/m_n(x_1,x_2) \mathrm{~for~} n \not = 1, 2$$
where $p_n$ is the weighted sum over all perfect matchings in $G_n$
given in definition \ref{14graphs}, with weighting as in the
preceding paragraph; and $m_n =x_1^{sq(n)}x_2^{oct(n)}$ with $sq(n)$
given by (\ref{sqn}) and $oct(n)$ given by (\ref{octn}).
\end{Thm}

In the above table, we see that Theorem \ref{14thm} is true for
small values of $n$, thus it suffices to prove that $p_n(x_1,x_2) /
x_1^{sq(n)}x_2^{oct(n)}$ satisfies the same periodic quadratic
recurrences as Laurent polynomials $x_n$.
%
By the definition of $m_n(x_1,x_2)$, it suffices to verify the
following two recurrences:
\begin{align}
\label{firststep} p_{2n+1}p_{2n+3} &= p_{2n+2} + x_1^{|4n+2|-2}x_2^{|2n|-1} \\
\label{secondstep} p_{2n}p_{2n+2} &= p_{2n+1}^4 +
x_1^{|8n|-4}x_2^{|4n-2|-2}.
\end{align}

\subsection{Proof of the first recurrence} \label{sectfirst}
We use a decomposition of superimposed graphs to prove equality
(\ref{firststep}). Unlike Kuo's technique of graphical condensation,
we will not use a central graph containing multi-matchings, but will
instead use superpositions that only
overlap on one edge, as in the bijective proof of Lemma
\ref{recur2}.
First, let $G_{2n+1}$ be defined as in the previous subsection, and
let $H_{2n + 3}$ (resp. $H_{2n-1}$) be constructed by taking graph
$G_{2n + 3}$ (resp. $G_{2n-1}$), reflecting it horizontally, and
then rotating the leftmost square upwards.  For convenience of
notation we will henceforth let $M = 2n+1$ so that we can abbreviate
these two cases as $H_{M\pm 2}$ (Throughout this section we choose
the sign of $H_{M\pm 2}$, $G_{M\pm 2}$ and $G_{M\pm 1}$ by using
$H_{2n+3}$, $G_{2n+3}$, and $G_{2n+2}$ if $n \geq 1$ and $H_{2n-1}$,
$G_{2n-1}$ and $G_{2n}$ if $n \leq -1$.) This reflection and
rotation will not change the number (or sum of the weights) of
perfect matchings. Thus the sum of the weights of perfect matchings
in the graph $H_{M \pm 2}$ also equals $p_{M \pm 2}$. Also, we will
let $K_2$ denote the graph
consisting of two vertices and a single edge connecting them.

The graph $G_{M \pm 1}$ can be decomposed as the union of the graphs
$G_{M}$, $H_{M \pm 2}$ and $K_2$ where $G_{M}$ and $H_{M \pm 2}$ are
joined together on an overlapping edge (the rightmost edge of
$G_{M}$ and the leftmost edge of $H_{M \pm 2}$). The graph $K_2$ is
joined to these graphs so that it connects to the bottom-right
(bottom-left) vertex of the rightmost octagon of $G_{M}$ and the
top-right (top-left) vertex of the leftmost octagon of $H_{2n+3}$
($H_{2n-1}$). See the picture below for an example with $G_{10}$.
The blue arc in the middle represents $K_2$. We will denote this
decomposition as
\begin{eqnarray*}
G_{M \pm 1} &=& G_{M} \cup H_{M \pm 2} \cup K_2.
\end{eqnarray*}

\vspace{1em}\begin{center}
\includegraphics[width = 3.2in , height = 0.8in]{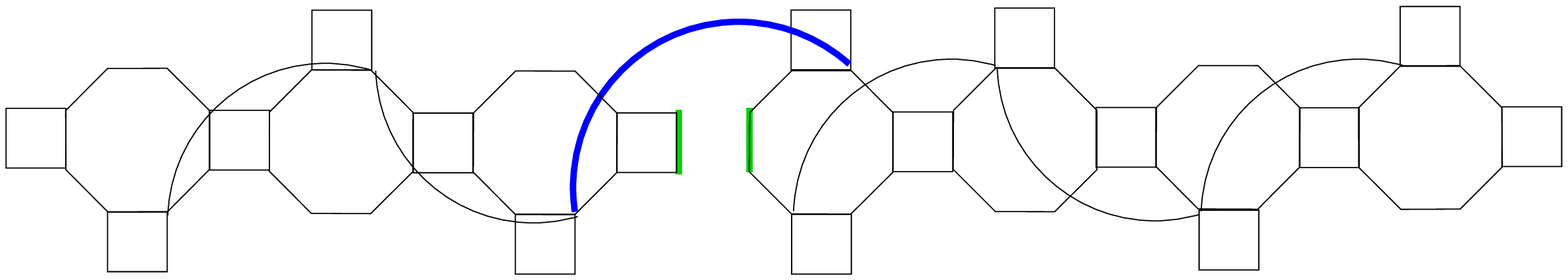}
\end{center}

As in subsection \ref{bij}, we let $G \sqcup H$ be the graph formed
by the disjoint union of graph $G$ and graph $H$. A perfect matching
of $G_{M}$ and a perfect matching of $H_{M \pm 2}$ will meet at the
edge of incidence in one of four ways (verticals meeting,
horizontals meeting vertical, verticals meeting diagonals,
or horizontals meeting diagonals).

In three of the cases, edges of weight 1 are utilized, and we can
bijectively associate a perfect matching of $G_{M} ~\sqcup~ H_{M \pm
2}$ to a perfect matching of $G_{M \pm 1}$ by removing an edge of
weight 1 on the overlap; though it is impossible to map to a perfect
matching of $G_{M \pm 1}$ that uses the edge of $K_2$ in this way.
This bijection is analogous to the one discussed in subsection
\ref{bij}. Thus we have a weight-preserving bijection between
\noindent $\{$perfect matchings of $G_{M \pm 1} ~- ~ K_2~\}$ and the
set $\{$perfect matchings of $G_{M} ~\sqcup~ H_{M \pm 2}\} ~-~
\{$pairs with nontrivial incidence (horizontals meeting diagonals)
$\}$ where by abuse of notation we here and henceforth let $G - K_2$
refer to the subgraph of $G$ with $K_2$'s edge deleted (without
deleting any vertices). Thus proving $p_{2n+1}p_{2n+3} - p_{2n+2} =
p_{M}p_{M+2} - p_{M+1} = x_1^{|4n+2|-2}x_2^{|2n|-1}$ reduces to
proving the following claim.

\begin{Prop} \label{heartoffirst}
The sum of the weights of all perfect matchings of $G_{M \pm 1}$
that contains $K_2$
is $x_1^{|4n+2|-2}x_2^{|2n|-1}$ less than the sum of the weights of
all perfect matchings of $G_{M} ~\sqcup~ H_{M \pm 2}$ that have
nontrivial incidence.
\end{Prop}

\noindent Before giving the proof of this Proposition we introduce a
new family of graphs that will allow us to write out several
steps of this proof more elegantly.  For $n \geq 1$, we let
$\tilde{G}_{2n+1}$ be the graph obtained from $G_{2n+1}$ by
\emph{deleting} the outer square on the extreme right. If $n \leq
-1$, we let $\tilde{G}_{2n+1}$ be the graph obtained from $G_{2n+1}$
by \emph{adjoining} an outer square on the extreme right.  For
example, $\tilde G_{9}$ and $\tilde G_{-7}$ are shown below. We let
$\tilde{p}_{2n+1}$ be the sum of the weights of perfect matchings in
$\tilde{G}_{2n+1}$.  Notice that this construction creates a
reciprocity such that graphs $\tilde{G}_{-M} = \tilde{G}_{-2n-1}$
and $\tilde{G}_{M+4} = \tilde{G}_{2n+5}$ are isomorphic.

\vspace{2em}
\includegraphics[width = 4in , height = 2in]{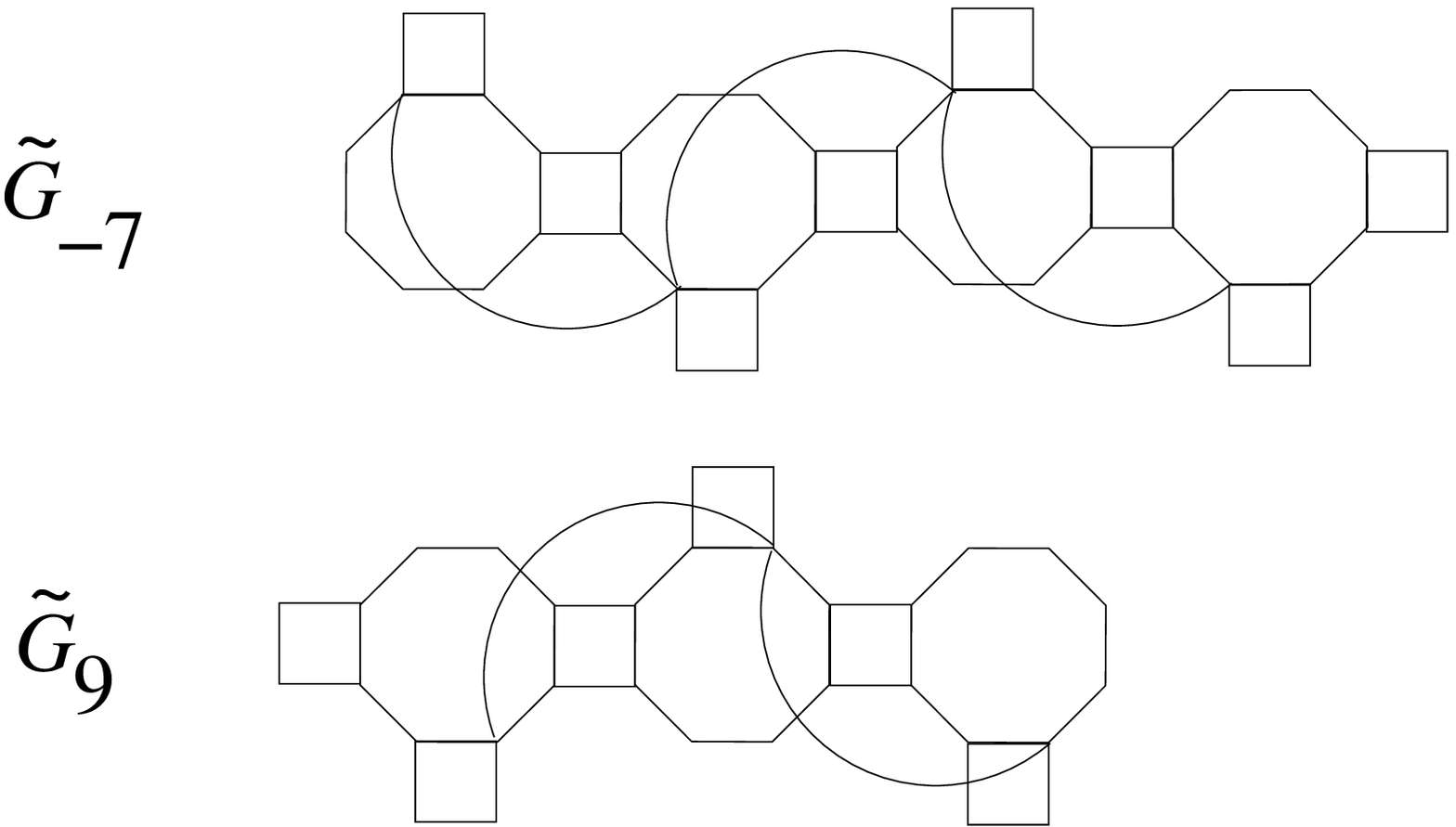}

\vspace{0.5em} \noindent The polynomials $p_{2n+1}$ and
$\tilde{p}_{2n+1}$ are related in a very simple way.

\begin{Lem} \label{Lem1}
\begin{align}
\label{tilderecur} p_{2n+1} &= (x_2+1)\tilde{p}_{2n+1} -
x_1^4x_2\tilde{p}_{2n-1} &\mathrm{for~} n &\geq 2, \\
\label{negtilderecur} p_{2n+1}  &= (x_1^4+x_2+1) \tilde{p}_{2n+3} -
x_1^4x_2^2\tilde{p}_{2n+5} &\mathrm{for~} n &\leq -2.
\end{align}
\end{Lem}

\begin{proof} The proof of Lemma \ref{Lem1} follows the same logic as the
inclusion-exclusion argument of section \ref{bij} that proved Lemma
\ref{recur2}. In this
case, we use the fact that we can construct $G_{2n+1}$ by adjoining
the graph $G_3$ (resp. $G_{-1}$) to $\tilde{G}_{2n+1}$ (resp.
$\tilde{G}_{2n+3}$) on the right for $n \geq 2$ (resp. $n \leq -2$).
It is clear that $p_3 = x_2+1$ (resp. $p_{-1} = x_1^4+x_2+1$) and
the only perfect matchings we must exclude are those that contain a
pair of diagonals
meeting a pair of horizontals.  For the $n \leq -2$ case, we must
also add in those perfect matchings that use the rightmost arc of
$G_{2n+1}$.
\end{proof}

\begin{proof} (Prop. \ref{heartoffirst})
By analyzing how the inclusion of certain key edges in a perfect
matching dictates how the rest of the perfect matching must look, we
arrive at the expressions
\begin{align}
\label{nontriv} x_1^4x_2(x_2+1)p_{2n-1}\tilde{p}_{2n+1} +
x_1^8x_2^2\tilde{p}_{2n-3}\tilde{p}_{2n+1} ~~~& \mbox{for $n \geq 2$} \\
\label{nontrivneg}  x_1^4x_2 p_{2n+3}\tilde{p}_{2n+1} +
x_1^4x_2^2\tilde{p}_{2n+5}\tilde{p}_{2n+1} ~~~& \mbox{for $n \leq
-2$}
\end{align}

\noindent for the sum of the weights of all perfect matchings of
$G_{M} \sqcup H_{M\pm 2}$ with nontrivial incidence (horizontals
meeting diagonals), and expressions
\begin{align}
\label{midarc}  x_1^4x_2(x_2+1)\tilde{p}_{2n-1}p_{2n+1} +
x_1^8x_2^2\tilde{p}_{2n-1}^2 ~~~& \mbox{for $n \geq 2$} \\
\label{midarcneg} x_1^4x_2   \tilde{p}_{2n+3}p_{2n+1} +
x_1^4x_2^2\tilde{p}_{2n+3}^2 ~~~& \mbox{for $n \leq -2$}
\end{align}

\noindent for the sum of the weights of all perfect matchings of
$G_{M\pm 1}$ using the $K_2$'s edge.

To prove Proposition \ref{heartoffirst}, it suffices to prove
$(\ref{nontriv}) = (\ref{midarc}) + x_1^{4n}x_2^{2n-1}$ and
$(\ref{nontrivneg}) = (\ref{midarcneg}) + x_1^{-4n-4}x_2^{-2n-1}$.  To
do so, we prove equalities

\begin{align}
\label{eq1} x_1^4x_2(x_2+1)(p_{2n-1}\tilde{p}_{2n+1} -
p_{2n+1}\tilde{p}_{2n-1})
&= x_1^{4n}x_2^{2n-2} + x_1^{4n}x_2^{2n-1}   \\
\label{eq2} x_1^8x_2^2(\tilde{p}_{2n-1}^2 -
\tilde{p}_{2n+1}\tilde{p}_{2n-3}) &= x_1^{4n}x_2^{2n-2}
\end{align}

\noindent  for $n \geq 2$ and equalities

\begin{align}
\label{eq3} x_1^4x_2(p_{2n+3}\tilde{p}_{2n+1} -
p_{2n+1}\tilde{p}_{2n+3})
&= x_1^{-4n}x_2^{-2n} + x_1^{-4n}x_2^{-2n+1}   \\
\label{eq4} x_1^4x_2^2(\tilde{p}_{2n+3}^2 -
\tilde{p}_{2n+1}\tilde{p}_{2n+5}) &= x_1^{-4n}x_2^{-2n}
\end{align}

\noindent for $n \leq -2$ and then subtract $(\ref{eq1}) -
(\ref{eq2})$ and $(\ref{eq3}) - (\ref{eq4})$. After shifting indices
and dividing both sides of equations (\ref{eq1}) through (\ref{eq4})
to normalize, we obtain that it suffices to prove

\begin{Lem} \label{midlem}

\begin{alignat}{2}
\label{mixed} p_{2n - 1}\tilde{p}_{2n+1} - p_{2n+1}\tilde{p}_{2n -
1}
&= x_1^{4n-4}x_2^{2n-3} & &\mbox{for $n \geq 2$} \\
\label{mixedneg} &= -x_1^{-4n}x_2^{-2n+1}(x_2+1)~~~~~ & &\mbox{for $n \leq -2$}   \\
\label{tildes} \tilde{p}_{2n+1}^2 - \tilde{p}_{2n-1}\tilde{p}_{2n +
3}
&= x_1^{4n-4}x_2^{2n-2} & &\mbox{for $n \geq 2$} \\
\label{tildesneg}  &= x_1^{-4n}x_2^{-2n} & &\mbox{for $n \leq -2$}.
\end{alignat}
\end{Lem}

We prove equations (\ref{tildes}) and (\ref{tildesneg})
simultaneously since $\tilde{p}_{2n+1} = \tilde{p}_{-2n+3}$. We
prove (\ref{tildes}) by making superimposed graphs involving
$\tilde{G}_{2n-1}$ and $\tilde{G}_{2n+3}$ and comparing it to the
superimposed graph of $\tilde{G}_{2n+1}$ with itself. Perhaps Kuo's
technique could be adapted to prove (\ref{tildes}), but instead we
consider a superposition overlapping over one edge, as we did
earlier in this proof.  Our superimposed graph thus resembles
$\tilde{G}_{4n-1}$, with a double edge somewhere in the middle and a
missing arc.  Analogous to the analysis that allowed us to reduce
from recurrence (\ref{firststep}) to Proposition \ref{heartoffirst},
we reduce our attention to the cases where gluing together
$\tilde{G}_{2n-1}$ and $\tilde{G}_{2n+3}$ and decomposing back into
two copies of $\tilde{G}_{2n+1}$ would not be allowed (or vice
versa). This entails focusing on cases where horizontals meet
diagonals at the double edge, or the arc appearing exclusively in
that decomposition (and not the other) appears in the matching.

After accounting for the possible perfect matchings, we find the
following two expressions
\begin{align}
\label{tild1} \tilde{p}_{2n-1}^2x_1^4x_2(x_2+1)
&+ \tilde{p}_{2n-3}\tilde{p}_{2n+1}x_1^4x_2   \\
\label{tild2} \tilde{p}_{2n-3}\tilde{p}_{2n+1}x_1^4x_2(x_2+1) &+
\tilde{p}_{2n-1}^2 x_1^4x_2
\end{align}

\noindent which represent the sum of the weights of nontrivial
perfect matchings in the superpositions of $\tilde{G}_{2n+1}$ and
$\tilde{G}_{2n+1}$ (resp. $\tilde{G}_{2n-1}$ and
$\tilde{G}_{2n+3}$). Taking the difference of these two expressions,
we find that
$$\tilde{p}_{2n+1}^2 - \tilde{p}_{2n-1}\tilde{p}_{2n+3}
= x_1^4x_2^2(\tilde{p}_{2n-1}^2 -
\tilde{p}_{2n-3}\tilde{p}_{2n+1}).$$ So after a simple check of the
base case ($\tilde{p}_5^2 - \tilde{p}_3\tilde{p}_7 = x_1^4x_2^2$) we
get equation (\ref{tildes}) by induction.

We easily derive (\ref{mixed}) from equations (\ref{tilderecur}) and
(\ref{tildes}):

\begin{align*}
p_{2n-1}\tilde{p}_{2n+1} - p_{2n+1}\tilde{p}_{2n-1} &=\\
( (x_2+1)\tilde{p}_{2n-1} -
x_1^4x_2\tilde{p}_{2n-3})\tilde{p}_{2n+1} - (
(x_2+1)\tilde{p}_{2n+1} - x_1^4x_2\tilde{p}_{2n-1})\tilde{p}_{2n-1} &= \\
x_1^4x_2(\tilde{p}_{2n-1}^2 - \tilde{p}_{2n-3}\tilde{p}_{2n+1}) &=
x_1^{4n-4}x_2^{2n-3}.
\end{align*}

\noindent We can also derive (\ref{mixedneg}) from (\ref{tildesneg})
but we first need to prove the following Lemma:

\begin{Lem} \label{keystep}
\rm \begin{alignat}{2} \label{keyst}
\tilde{p}_{2n+1}\tilde{p}_{2n-1} - \tilde{p}_{2n-3}\tilde{p}_{2n+3}
&= x_1^{4n-8}x_2^{2n-4}(x_1^4+(x_2+1)^2) & ~~~&\mbox{for $n \geq 2$}, \\
\label{keystneg} \tilde{p}_{2n+1}\tilde{p}_{2n-1} -
\tilde{p}_{2n-3}\tilde{p}_{2n+3} &=
x_1^{-4n}x_2^{-2n}(x_1^4+(x_2+1)^2) & ~~~ &\mbox{for $n \leq -2$}.
\end{alignat}
\end{Lem}

\begin{proof}We prove this Lemma by using the the same technique that we
used to prove equation (\ref{tildes}).  Analogously, they can be
proven simultaneously by proving (\ref{keyst}) because of the
reciprocity $\tilde{p}_{2n+1} = \tilde{p}_{-2n+3}$.  In this case,
the superimposed graph resembles $\tilde{G}_{4n-1}$ and by
considering nontrivial perfect matchings, we obtain that
$$\tilde{p}_{2n+1}\tilde{p}_{2n-1} -
\tilde{p}_{2n-3}\tilde{p}_{2n+3} =
x_1^4x_2^2(\tilde{p}_{2n-1}\tilde{p}_{2n-3} -
\tilde{p}_{2n-5}\tilde{p}_{2n+1}).$$  Since

$$\tilde{p}_{5}\tilde{p}_{7} -
\tilde{p}_{3}\tilde{p}_{9} = x_1^4x_2^2(x_1^4+(x_2+1)^2),$$

\noindent we have the desired result by induction.\end{proof}

We thus can verify (\ref{mixedneg}) algebraically by using
(\ref{negtilderecur}), (\ref{tildesneg}), and (\ref{keystneg}).
Since the proof of equations (\ref{mixed}) through (\ref{tildesneg})
was sufficient, thus recurrence (\ref{firststep}) is proven for $n
\not = 0$ or $-1$.
\end{proof}

\subsection{Proof of the second recurrence} \label{sectsecond}
We now prove the recurrence (\ref{secondstep}) via the following two
observations.

\begin{Lem} \label{Prop2} \rm
\begin{alignat}{2} \label{Lem2}
p_{2n-1} p_{2n+3} - p_{2n+1}^2 &=
x_1^{4n-4}x_2^{2n-3}(x_1^4+(x_2+1)^2) &~~~~ &\mbox{for $n \geq 2$} \\
\label{negLem2} &= x_1^{-4n-4}x_2^{-2n-1}(x_1^4+(x_2+1)^2) &~~~~ &\mbox{for
$n \leq -1$}.
\end{alignat}
\end{Lem}

\begin{Lem} \label{incexc} \rm
\begin{alignat}{2}
\label{Lem3} p_{2n+3} &= (x_1^4+(x_2+1)^2 )p_{2n+1} -
x_1^4x_2^2 p_{2n-1} &~~~~ &\mbox{for $n \geq 2$}, \\
\label{negLem3} p_{2n-1} &= (x_1^4+(x_2+1)^2 )p_{2n+1}
-x_1^4x_2^2p_{2n+3} &~~~~ &\mbox{for $n \leq -1$}.
\end{alignat}
\end{Lem}

\begin{proof} (Lemma \ref{Prop2}) We easily derive (\ref{Lem2}) from Lemmas \ref{Lem1},
\ref{midlem} and \ref{keystep} by the following derivation,
{\footnotesize
\begin{align*}
p_{2n-1} &p_{2n+3} - p_{2n+1}^2 = \\
&((x_2+1)\tilde{p}_{2n-1} -
x_1^4x_2\tilde{p}_{2n-3})((x_2+1)\tilde{p}_{2n+3} -
x_1^4x_2\tilde{p}_{2n+1}) - ((x_2+1)\tilde{p}_{2n+1} -
x_1^4x_2\tilde{p}_{2n-1})^2
= \\
&(x_2+1)^2(\tilde{p}_{2n-1}\tilde{p}_{2n+3} - \tilde{p}_{2n+1}^2)
~~+ x_1^8x_2^2(\tilde{p}_{2n-3}\tilde{p}_{2n+1} -
\tilde{p}_{2n-1}^2) ~~-
x_1^4x_2(x_2+1)(\tilde{p}_{2n-3}\tilde{p}_{2n+3} -
\tilde{p}_{2n+1}\tilde{p}_{2n-1}) = \\
&-(x_2+1)^2(x_1^{4n-4}x_2^{2n-2})-x_1^8x_2^2(x_1^{4n-8}x_2^{2n-4}
+x_1^4x_2(x_2+1)(x_1^4+(x_2+1)^2)x_1^{4n-8}x_2^{2n-4}) =\\
& x_1^{4n-4}x_2^{2n-3}(x_1^4+(x_2+1)^2).
\end{align*}}
\noindent The proof of (\ref{negLem2}) is similar and thus Lemma
\ref{Prop2} is proved.\end{proof}

\begin{proof} (Lemma \ref{incexc}) Like Lemma \ref{Lem1}, Lemma \ref{incexc} can also be proven
using an inclusion-exclusion argument.  This one relies on the fact
that $G_{2n+3}$ is inductively built from $G_{2n+1}$ by adjoining an
octagon, an arc, and two squares.  Graph $\tilde{G}_{5} \sqcup
G_{2n+1}$ has $(x_1^4+(x_2+1)^2)p_{2n+1}$ as the sum of the weight
of all its perfect matchings. Most perfect matchings of
$\tilde{G}_{5} \sqcup G_{2n+1}$ map to a perfect matching of
$G_{2n+3}$ (resp. $G_{2n-1}$) with the same weight.  The only
perfect matchings that do not participate in the bijection are those
that contain a pair of diagonals meeting
a pair of horizontals. The sum of the weights of all such perfect
matchings is $x_1^4x_2(x_2+1)p_{2n-1}$ (resp.
$x_1^4x_2(x_2+1)p_{2n+3}$). However, we have neglected the perfect
matchings of $G_{2n+3}$ (resp. $G_{2n-1}$) that use the one arc not
appearing in $\tilde{G}_{5}
\sqcup G_{2n+1}$. Correcting for this we add back $x_1^4x_2p_{2n-1}$
(resp. $x_1^4x_2p_{2n+3}$). After these subtractions and additions
we do indeed obtain equation (\ref{Lem3}) for $n \geq 2$ (resp.
(\ref{negLem3}) for $n \leq -1$).\end{proof}

We now are ready to prove (\ref{secondstep}). Using the first
recurrence, (\ref{firststep}), we can rewrite the lefthand side of
(\ref{secondstep}) as

$$(p_{2n-1} p_{2n+1} - x_1^{|4n-2|-2}x_2^{|2n-2|-1})(p_{2n+1} p_{2n+3}
- x_1^{|4n+2|-2}x_2^{|2n|-1})$$

\noindent which reduces to

$$p_{2n+1}^2 (p_{2n-1}p_{2n+3})
- p_{2n+1}x_1^{|4n-2|-2}x_2^{|2n-2|-1} (p_{2n+3} + x_1^4x_2^2
p_{2n-1}) + x_1^{|8n|-4} x_2^{|4n-2|-2}.$$

\noindent Using Lemma \ref{Prop2} and Lemma \ref{incexc}, this
equation simplifies to $p_{2n+1}^4 + x_1^{|8n|-4}x_2^{|4n-2|-2}$.
Thus the recurrence (\ref{secondstep}) is proven for $n \not =  0$
or $1$.  We have thus proven Theorem \ref{14thm}
.

\subsection{A combinatorial interpretation for the semicanonical basis}
\label{14semi14}
It was shown in \cite{Sherman} that a canonical basis for the
positive cone consists of cluster monomials, that is monomials of
the form $x_n^px_{n+1}^q$, as well as one additional sequence of
elements, in the affine ($(2,2)$ or $(1,4)$) case.
One can think of this extraneous sequence as corresponding to the
imaginary roots of the Kac Moody algebra, which are of the form
$n\delta$ where $\delta = \alpha_1 + \alpha_2$ in the $(2,2)$ case
and $2\alpha_1+\alpha_2$ in the $(1,4)$ case.

As described in \cite{Sherman}, this sequence completes the
canonical basis, and is closely related to the sequence of $s_n$'s
which completes the \emph{semicanonical} basis. The $s_n$'s are
defined as the normalized Chebyshev polynomials of the second kind
in variable $z_1=s_1 = (x_1^4 + (x_2+1)^2)/x_1^2x_2$ just as they
were in the $(2,2)$-case \cite{Sherman}.  Using this definition, we
obtain a combinatorial interpretation for the $s_n$'s in the
$(1,4)$-case just as we did in the $(2,2)$-case.
In both cases, the non-cluster monomial elements of the
semicanonical basis were discovered as an auxiliary sequence in the
proof of $x_n$'s combinatorial interpretation.




\begin{Thm} In the $(1,4)$-case, the Laurent polynomials
$s_n(x_1,x_2)$, defined as $\overline{S}_n(z_1)$, are precisely
$\tilde{p}_{2n+3}(x_1,x_2)/\tilde{m}_{2n+3}(x_1,x_2)$ where the
$\tilde{p}_{2n+3}$'s are defined in section \ref{sectfirst} (between
Proposition \ref{heartoffirst} and Lemma \ref{Lem1}) and
$\tilde{m}_{2n+3}(x_1,x_2) = x_1^{2n}x_2^{n}.$
\end{Thm}

Therefore the $s_n$'s have a combinatorial interpretation in terms
of the graphs $\tilde{G}_{2n+1}$ of section \ref{sectfirst}.

\begin{proof}
Our method of proof is analogous to our proof of Lemmas \ref{recur2}
and \ref{Lem1}.  We note that a perfect matching of
$\tilde{G}_{2n+3}$ can be decomposed into a perfect matching of
$\tilde{G}_{2n+1}$ and $\tilde{G}_5$ (with the graph $\tilde{G}_5$
on the righthand side) or will utilize the rightmost arc.  However,
there is not a bijection between perfect matchings avoiding the
rightmost arc and perfect matchings of $\tilde{G}_{2n+1} \sqcup
\tilde{G}_5$ since we must exclude those matchings of
$\tilde{G}_{2n+1} \sqcup \tilde{G}_5$ that 
use the two rightmost horizontal edges of $\tilde{G}_{2n+1}$ and the
two diagonal edges of $\tilde{G}_5$.  In conclusion, we obtain

\begin{eqnarray}
\tilde{p}_{2n+3} &=& \tilde{p}_{2n+1}\tilde{p}_5 +
x_1^4x_2\tilde{p}_{2n-1} - x_1^4x_2(x_2+1)\tilde{p}_{2n-1} \\
\label{semieqn} &=& \tilde{p}_{2n+1}\tilde{p}_5 -
x_1^4x_2^2\tilde{p}_{2n-1}.
\end{eqnarray}
After dividing by $\tilde{m}_{2n+3}$, $\tilde{m}_{2n+1}$,
$\tilde{m}_{2n-1}$, and $\tilde{m}_{5}$ accordingly, equation
(\ref{semieqn}) reduces to

\begin{align} {\tilde{p}_{2n+3}\over \tilde{m}_{2n+3}} =
{\tilde{p}_{5}\over \tilde{m}_{5}} {\tilde{p}_{2n+1}\over
\tilde{m}_{2n+1}} - {\tilde{p}_{2n-1}\over \tilde{m}_{2n-1}}
\end{align}
and thus the $\tilde{p}_{2n+3}/\tilde{m}_{2n+3}$'s satisfy the same
recurrence as the normalized Chebyshev polynomials of the second
kind \cite[Equation (3.2)]{Caldero} or \cite[Equation (2)]{NewNote};
and thus the same recurrence as the $s_n$'s.
\end{proof}

\section{Comments and open problems} \label{openprobs}

In both the (2,2)- and (1,4)- cases, we have now shown that the
sequence of Laurent polynomials $\{x_n\}$ defined by the
appropriate recurrences have numerators with positive coefficients.
A deeper combinatorial understanding of the $(1,4)$ case might, in
combination with what is already known about the other cases of
$(b,c)$ with $bc \leq 4$, give us important clues into how one might
construct suitable graphs $G_n$ for the cases $bc>4$.

%
%
%
Another direction follows from the combinatorial interpretation of
semicanonical basis elements in the $(1,4)$-case.  This analysis
motivates the question of whether or not one can find a family of
graphs such that the sum of the weights of their perfect matchings
are precisely the numerators of the $z_n$'s.  This would give
combinatorial interpretation to not only the semicanonical basis,
but also the canonical basis.
%

\vspace{2em} \noindent \bf Acknowledgements. \rm  We thank Sergey
Fomin for many useful conversations, Andrei Zelevinsky for his
numerous suggestions and invaluable advice, and an anonymous referee for
helpful comments.
The computer programs {\tt graph.tcl} (created by Matt Blumb) and
{\tt Maple} were the source of many experiments that were invaluable
for finding the
patterns proved in this research.  We would also like to acknowledge
the support of the University of Wisconsin, Harvard University, the
NSF, and the NSA for their support of this project through the VIGRE
grant for the REACH program.

\end{document}